\documentclass[11pt]{article}
\usepackage[centertags]{amsmath}
\usepackage{amsmath}
\usepackage{amssymb}
\usepackage{amsthm}

\usepackage{color}

\numberwithin{equation}{section}

\newcommand{\1}{1\!\!1}
\newtheorem{definition}{Definition}[section]
\newtheorem{theorem}{Theorem}[section]
\newtheorem{proposition}{Proposition}[section]
\newtheorem{lemma}{Lemma}[section]
\newtheorem{remark}{Remark}[section]
\newtheorem{corollary}{Corollary}[section]
\newtheorem{example}{Example}[section]

\def \R{\mathbb{R}}

\def \E{\mathbb{E}}

\def \L{\mathbb{L}}

\def \bf{\textbf}
\def \it{\textit}

\def \bop {\noindent\textbf{Proof }}
\def \eop {\hbox{}\nobreak\hfill
	\vrule width 2mm height 2mm depth 0mm
	\par \goodbreak \smallskip}

\def \bop {\noindent\textbf{Proof }}
\def \eop {\hbox{}\nobreak\hfill
	\vrule width 2mm height 2mm depth 0mm \par \goodbreak \smallskip}

\usepackage{amssymb}
\usepackage{amsthm}
\usepackage{amsmath}
\usepackage{latexsym}
\usepackage[T1]{fontenc}
\usepackage[latin1]{inputenc}

\def \R{\mathbb{R}}

\def \E{\mathbb{E}}

\def \L{\mathbb{L}}

\def \bf{\textbf}
\def \it{\textit}

\def \bop {\noindent\textbf{Proof.}}
\def \eop {\hbox{}\nobreak\hfill
	\vrule width 2mm height 2mm depth 0mm
	\par \goodbreak \smallskip}

\def \eop {\hbox{}\nobreak\hfill \vrule width 2.0mm height 1.8mm depth 0mm
	\par \goodbreak \smallskip}

\numberwithin{equation}{section}

\def \R{\mathbb{R}}

\def \E{\mathbb{E}}

\def \L{\mathbb{L}}

\def \bf{\textbf}
\def \it{\textit}

\def \bop {\noindent\textbf{Proof }}
\def \eop {\hbox{}\nobreak\hfill
	\vrule width 2mm height 2mm depth 0mm
	\par \goodbreak \smallskip}

\begin{document}

	
	\title{      Solving Unbounded Quadratic BSDEs \\  
		by a  Domination  method  
		\thanks{Partially supported by
			PHC Toubkal/18/59.} }
	\author{Khaled Bahlali
		\\ \\
		Université de Toulon,
		IMATH,
		EA 2134, \\
		83957 La
		Garde
		Cedex, France. }
	\date{}
	\maketitle

	\noindent {\bf{Abstract}} \ We introduce a domination argument which asserts that:  if we can dominate the parameters  of a quadratic backward stochastic differential equation (QBSDE) with continuous generator from above and from below by those of two BSDEs having ordered solutions, then also the original  QBSDE admits at least one solution. This result is presented in a general framework: we do not impose any integrability condition on none of the terminal data of the three involved BSDEs,  we do not require any constraint on the growth nor continuity of the two dominating generators.
	As a consequence, we establish the existence of a maximal and a minimal  solution to   BSDEs whose coefficient $H$ is continuous and satisfies \
	 $\vert H(t,y,z)\vert\leq \alpha_t
	+ \beta_t\vert y\vert + \theta_t\vert z\vert + f(|y|)\vert z\vert^{2}$, where $\alpha_t$, $\beta_t$, $\theta_t$  are positive processes and the function $f$ is positive, continuous and increasing (or even  only positive and locally bounded) on $\R$. This is done with unbounded terminal value. We cover the classical QBSDEs where the function $f$ is constant (\cite{BEK}, \cite{BH1},  \cite{Kob}, \cite{LSM2}) and when $f(y) = y^p$ (\cite{fd}) and  also the cases  where the generator has super linear growth such as \ $y|z|$,   $e^{|y|^k} |z|^p$, $e^{e^{|y|}} |z|^2$,   ($ k \geq 0$, $0\leq p < 2$) and so on.  In contrast to the works \cite{BEK, BH1, fd, Kob, LSM2}, we get the existence of a  a maximal and a minimal solution and we cover the BSDEs with at most linear growth  (take $f=0$). In particular, we cover and extend the results of \cite{ht} and \cite{LSM1}. Furthermore, we establish the existence and uniqueness of solutions  to BSDEs driven by $f(y)|z|^2$ when $f$ is merely locally integrable on $\R$.

	\vskip 0.2cm
	\noindent \text{AMS 2000 Classification subjects:} 60H10, 60H20, 60H30, 91G10
	
	\vskip 0.2cm \noindent {${Keywords}:$}  Superlinear backward stochastic
	differential equations, Quadratic backward stochastic
	differential equations, unbounded solutions.
	maximal and a minimal

	\medskip
	
	\section{Introduction}

Let $(W_{t})_{0\leq t\leq T}$ be a $d$--dimensional Brownian motion defined on a complete probability space $(\Omega ,\mathcal{F}, \mathbb{P})$. We
denote by $(\mathcal{F}_{t})_{0\leq t\leq T}$ the natural filtration of $W$
augmented with $\mathbb{P}$-negligible sets. Let $H(t,\omega, y, z)$ be a
real valued $\mathcal{F}_{t}$-progressively measurable process defined on $
[0, \ T]\times\Omega\times\mathbb{R}\times\mathbb{R}^d$. Let $\xi$ be an $
\mathcal{F}_{T}$-measurable $\mathbb{R}$-valued random variable. Consider the BSDE
\begin{equation}\label{H}
Y_{t}=\xi +\int_{t}^{T}H(s,Y_{s},Z_{s})ds-\int_{t}^{T}Z_{s}dW_{s},\, \ \ 0\leq
t\leq T \tag{$eq(\xi ,H)$}
\end{equation}
$\xi$ is called the terminal value and $H$ is called the generator or the coefficient. A BSDE with data $(\xi, H)$ will be labeled  $eq(\xi, H)$ or BSDE $(\xi, H)$ or  BSDE$(\xi, H)$.

\begin{definition}
		\label{defsol} $(i)$ \ We say that $eq(\xi ,H)$ is quadratic if $H$ has at most a quadratic growth in its $z$-variable. 	
		
			\vskip  0.2cm
		$(ii)$ \ A solution to $eq(\xi , H)$ is a process
		$(Y,Z)$ which satisfies $eq(\xi , H)$ on $ [0, \ T]$ and such that $Y$ is continuous, $ \int_0^T |Z_s|^2ds < \infty$ \ $a.s$ and \ $\int_0^T |H(s, Y_s, Z_s)|ds < \infty$ \ $a.s$.
		
		\vskip 0.2cm
		$(iii)$ \ A positive solution is a solution $(Y,Z)$ such that $Y_t \geq 0$. We Symmetrically  define a negative solution. A bounded solution is a solution $(Y,Z)$ such that $Y$ is bounded
	\end{definition}
$eq(\xi, H)$ is related to partial differential equations (PDEs), optimal stochastic control and mathematical finance (risque measure, utility maximization, etc.). The Quadratic BSDEs  were  studied in many papers, among them one can cite the works \cite{beo2013, beo2017, BEK,  BH1,  DHR1, deps, EH2011, fd, Kob, LSM2, Tevz}.

In this paper, we are concerned with the existence of solutions to BSDEs whose generator $H$ satisfies $\vert H(t,y,z)\vert\leq \alpha_t
	+ \beta_t\vert y\vert + \theta_t\vert z\vert + f(|y|)\vert z\vert^{2}$, where $\alpha_t$, $\beta_t$, $\theta_t$  are positive processes and the function $f$ is  positive on $\R_+$ and locally bounded but not globally integrable on $\R$. We are motivated by the fact that the BSDEs driven by $H(t, y, z) =  g(c_t,y) +  f(|y|)\vert z\vert^{2}$  appears in stochastic differential utility, see  \cite{deps}. This type of BSDEs are also related  to quadratic PDEs appearing in financial markets, see \cite{DJ}. Let us present another motivation : it has been recently shown in \cite{beo2013, beo2017} that the BSDEs driven by a generator $H$ satisfying $\vert H(t,y,z)\vert\leq \alpha
	+ \beta\vert y\vert + \theta\vert z\vert + f(|y|)\vert z\vert^{2}$ have solutions when $f$ is globally integrable and $\alpha$, $\beta$, $\theta$ are constant. However, these two works can not cover the classical BSDEs driven by  $|z|^2$, since they assume that $f$ is integrable on $\R$ . Thus, the following questions naturally arise:

\textit{1) Are there BSDEs whose generator $H$  satisfying  $\vert H(t,y,z)\vert\leq \alpha_t
	+ \beta_t\vert y\vert + \theta_t\vert z\vert + f(|y|)\vert z\vert^{2}$ that have solutions without assuming the global integrability of $f$ ?}

\textit{2) If yes, what integrability condition we should require on the terminal value  $\xi$ ?}

\vskip 0.2cm  The following example gives a positive answer to the first question. It moreover  shows  that neither global integrability of $f$ nor integrability of $\xi$ are necessary to the existence of solutions.

\vskip 0.2cm \noindent
\begin{example}  {\rm {Consider the BSDE
		\begin{equation}\label{exintro}
	Y_{t}=\xi +\int_{t}^{T}\1_{\R_+}(Y_s)|Z_s|^2ds-\int_{t}^{T}Z_{s}dW_{s},\, \ \ 0\leq
	t\leq T.
		\end{equation}
		where for a set $A$,  $\1_{A}$ denotes the indicator function of $A$.
		
   Equation \eqref{exintro} is not covered by the works \cite{BEK, BH1, deps, DHR1, EH2011, fd, Kob, LSM2, Tevz}, since $\1_{\R_+}$ is neither constant nor continuous. It is also not covered by the works \cite{beo2013, beo2017, deps}, since the function $\1_{\R_+}$ is not globally integrable.  Nevertheless, equation \eqref{exintro}  admits a solution \textit{without any integrability} condition on $\xi$.  Indeed,
 the function $u(y) := \frac12(e^{2y}- 1)\1_{\R_+}(y) + y\1_{\R_-^*}(y)$ belongs to the Sobolev space $W _{1,\,loc}^{2}(\R)$ and solves the differential equation $\frac12 u^{\prime\prime}(y) - \1_{\R_+}(y)u^{\prime}(y) = 0$   on $\R^*$. Therefore, using Itô-Krylov's formula for BSDEs (see \cite{beo2013, beo2017}) one can show that equation \eqref{exintro} has a solution  if and only if the following equation has a solution.
 \begin{equation}\label{exintro0}
	\bar Y_{t}=u(\xi )-\int_{t}^{T}\bar Z_{s}dW_{s},\, \ \ 0\leq
	t\leq T.
		\end{equation}
By Dudley's representation theorem \cite{Dud}, equation \eqref{exintro0} has a solution without any integrability condition on $u(\xi )$. Since $u$ is a one to one function from $\R$ onto $\R$,  we deduce that equation \eqref{exintro}  has a solution \textit{without any integrability of  $\xi$}. If moreover, $u(\xi )$ is integrable then equation \eqref{exintro} has a unique solution $(Y,Z)$ such that $u(Y)$ belongs to class $(D)$. Other examples will be presented later. }}
\end{example}

The aim of this paper is to introduce a domination method  which allows to solve    $eq(\xi, H)$ when $\xi$ is unbounded and $\vert H(t,y,z)\vert\leq \alpha_t
+ \beta_t\vert y\vert + \theta_t\vert z\vert + f(|y|)\vert z\vert^{2}$, where $\alpha_t$, $\beta_t$, $\theta_t$  are positive processes and the function $f$ is  positive on $\R_+$ and locally bounded but not globally integrable on $\R$. Our strategy is divided into five stages each of which has its own interest:  


In the first step $(i)$, we establish  existence, uniqueness and comparison of solutions to BSDEs driven by $f(y)|z|^2$ when $f$ is locally integrable, this covers the classical case where $f$ is constant. It will be shown that when $u_f(\R) = \R$, then the BSDE $(\xi, f(y)|z|^2)$  has a solution without any integrability condition on the terminal value, and when $u_f(\R) \neq \R$,  then the condition  $u_f(\xi)$ \textit{integrable} is necessary to the existence of solutions for $eq(\xi, f(y)|z|^2)$. We also show that the uniqueness as well as the comparison hold for   BSDE$(\xi, f(y)|z|^2)$ in the class of solutions such that $u_f(Y)$ belongs to class $(D)$.


\vskip 0.1cm
In order to explain the other steps, we precise some notation, definitions and assumptions.

\vskip 0.2cm \noindent
	 \textbf{Some notation.} For given real numbers $a$ and $b$, we set \ $a\wedge b:=\min (a,b)$, \ $
	a\vee b:=\max (a,b)$, \ $a^{-}:=\max (0,-a)$ \ and \ $a^{+}:=\max (0,a)$.

\vskip 0.15cm
For given positive processes $\alpha$ and $\beta$ we denote \ $\xi_{\alpha, \beta} := \left(|\xi|+ \int_0^T\alpha_s ds\right)e^{\int_0^T\beta_s ds}$.  	

We define $\xi_{\alpha, \beta}^+$
\  and \ 		
$\xi_{\alpha, \beta}^-$  in likewise manner.

\vskip 0.15cm
	For $p>0$, we denote by \ $\L _{loc}^{p}(\R)$ ($\L _{loc}^{p}$ in short) the
	space of (classes) of functions $u$ defined on $\mathbb{R}$ which are $p$-integrable on bounded set of $\R$.
	We
	also denote,
	
	\vskip 0.15cm
	$W _{p,\,loc}^{2}$ := the Sobolev space of (classes) of
	functions $u$ defined on $\mathbb{R}$ such that both $u$ and its generalized
	derivatives $u^{\prime }$ and $u^{\prime \prime }$ belong to $\mathbb{L}
	_{loc}^{p}(\mathbb{R})$.
	
	\vskip 0.15cm
	$\mathcal{C}$ \ := \ the space of continuous and $\mathcal{F}_{t}$
	--adapted processes.
	
	\vskip 0.15cm
	$\mathcal{S}^{p}$ := the space of continuous, $\mathcal{F}_{t}$
	--adapted processes $\varphi$ such that \ $\mathbb{E}\left(\sup_{0\leq
		t\leq T}|\varphi_{t}|^{p}\right) < \infty. $
	
	\vskip 0.15cm
	$\mathcal{L}^{2}$ := the space of $\mathcal{F}_{t}$ --adapted processes $
	\varphi$ satisfying $\int_{0}^{T}|\varphi_{s}|^{2}ds<+\infty \ \mathbb{P}
	\text{--a.s.} \label{inq2} $
	
	$\mathcal{M}^{p}$:= the space of $\mathcal{F}_{t}$--adapted processes $
	\varphi$ satisfying $\mathbb{E}\left[ \left(\int_{0}^{T}|\varphi_{s}|^{2}ds\right)^{\frac{p}{2}}\right]<+\infty. $

   \vskip 0.15cm
    We say that the process $ \varphi:=(\varphi_s)_{0\leq s\ T}$ \ belongs to class $(D)$ if
 $\sup_{0\leq\tau\leq T}|Y_{\tau}| < \infty$, where the supremum is taken over all stopping times $\tau$ such that  $\tau \leq T$.

\vskip 0.15cm $BMO$ \ is the space of uniformly integrable martingales $M$ satisfying $\sup_\tau||\E\big(|M_T-M_\tau| / {\cal F}_\tau\big)||_{\infty}<\infty$,
where the supremum is taken over all stopping times $\tau$.

 \vskip 0.15cm
 For a given locally integrable function $f$ defined on $\R$, we put
\begin{equation}\label{u}
u_f(y):\int_{0}^{y}\exp \left( 2\int_{0}^{x}f(r)dr\right) dx.
\end{equation}
The properties of the function $u_f$ are given in Lemma \ref {zvonkinM} of Apendix.

\vskip 0.2cm Consider the following assumptions.

\vskip 0.2cm \noindent	\textbf{(A1)} \hskip 0.3cm $H$ is continuous in $(y,z)$ for $a.e$ $(t,\omega)$ \ and satisfies,	
	\begin{equation*} \vert H(t,y,z)\vert\leq \alpha_t + \beta_t\vert y\vert  + f(|y|)\vert z\vert^{2}
\end{equation*}
	\textit{where $\alpha_t$, $\beta_t$ are some $(\mathcal{F}_t)$-adapted processes which are positive and $f$ is a real valued  function which is continuous, increasing and positive on $\mathbb{R}_+$}
$$
 \leqno{\textbf{(A2)}}  \qquad \qquad \qquad   \displaystyle u_f\left(\xi_{\alpha, \beta}\right) \  is \  integrable,
 $$
 where $u_f$ is defined by \eqref{u}.

\vskip 0.2cm
In step $(ii)$, we introduce  the domination argument  (Lemma \ref{d} below) which is an abstract result that gives the existence of solutions to BSDEs without integrability condition on the terminal value. It asserts that:  \textit{if  $(\xi_1, H_1)$ and $(\xi_2, H_2)$ are two BSDEs which respectively admit two solutions $(Y_1, Z_1)$ and $(Y_2, Z_2)$ such that   $(\xi_1, H_1, Y_1) \leq  (\xi_2, H_2, Y_2)$, then any Quadratic BSDE $(\xi, H)$ with continuous generator and satisfying  $(\xi_1, H_1) \leq (\xi, H) \leq  (\xi_2, H_2)$ has at least one solution $(Y,Z)$ such that $ Y_1 \leq Y \leq  Y_2$. Moreover, among all solutions which lie between $Y_1$ and $Y_2$, there are a maximal and a minimal solution}. The proof of this result do not need any a priori estimate nor approximation. It is based on the remarkable work  \cite{EH2013} on the existence of reflected QBSDEs without any integrability condition on the terminal value. Actually, we derive the existence of our QBSDE from a suitable reflected QBSDE.

In steps $(iii)$--$(iv)$, we use the domination   argument to show that when assumptions (A1), (A2) are satisfied then $eq(\xi, H)$  admits a \textit{maximal and a minimal  solution} satisfying some conditions which will be precised later.   We  cover the results obtained in \cite{BEK, BH1, fd, Kob} and also the cases   where the generator has super linear growth such as \ $y|z|$,  $e^{|y|}|z|^{p}$ ($0\leq p < 2$) $e^{|y|}|z|^{2}$ and more generally  the case  $|H(y,z)| \leq \alpha_t + \beta_t|y| + g(|y|)|z|^{p} + f(|y|)|z|^2 $ where  $f$ and $g$  are continuous, increasing and positive on $\R_+$, and $0\leq p < 2$. In all these cases, the existence of solutions is established with an unbounded terminal value.
Although the case where the function  $f$ is globally integrable  is not covered by assumption (A1), one can again use the domination argument to show that if $H$ satisfies (A1) with $f$ globally integrable and  $\xi_{\alpha, \beta}$ is  integrable,
then $eq(\xi, H)$ has a maximal and a minimal  solution; this extends the result of \cite{beo2013, beo2017} to the case where  $\alpha$, $\beta$ are processes and $\xi$ is merely in $\mathbb{L}^1$.

In step $(v)$, we establish the existence of solutions to $eq(\xi, H)$ under the the following assumptions:

\vskip 0.2cm \noindent	\textbf{(A3)} \hskip 0.3cm $H$ is continuous in $(y,z)$ for $a.e$ $(t,\omega)$ \ and satisfies,	
\begin{equation*} \vert H(t,y,z)\vert\leq \alpha_t + \beta_t\vert y\vert  + \theta_t\vert z\vert+ f(|y|)\vert z\vert^{2}
\end{equation*}
\textit{where $\alpha_t$, $\beta_t$ are some $(\mathcal{F}_t)$-adapted processes which are positive and $f$ is a real valued  function which is continuous, increasing and positive on $\mathbb{R}_+$}

$$
\leqno{\textbf{(A4)}}  \qquad \displaystyle  \sup_{\pi\in\sum} \E\left(\Gamma_{0, T}^{\pi}u_f(\xi_{\alpha, \beta})\right) :=\sup_{\pi\in\sum} \E\left({}e^{\int_0^T \theta_r \pi_rdW_r -\frac12 \int_0^T
	\theta_r^2 |\pi_r|^2dr}u_f(\xi_{\alpha, \beta})\right)  < +\infty
$$

where   \
$ \sum := \left\{  \pi\in {\cal L}^{2},\, |\pi|\in  \{0, 1\},  \  a.e. \ {\rm{and}} \  \displaystyle \  ess\sup_{\omega} \int_0^T
\theta_r^2 |\pi_r|^2dr < +\infty\right\}$. 	

\vskip 0.15cm\noindent In this case, we use again the domination argument to reduce the solvability of  $eq(\xi, H)$ to that of  $eq\big(u_f(\xi_{\alpha, \beta}\big),\, \theta_t|z|)$ and then to deduce the existence of solutions to $eq(\xi, H)$. It should be noted that the works \cite{beo2013, beo2017, BEK, BH1, deps, DHR1, fd, Kob, LSM2, Tevz} consider only the case $\theta_t = 0$. To the best of our knowlege, the case $\theta_t \neq 0$ is considered only  in \cite{beo2017, EH2011} and in the present paper. We emphasize that assumptions (A1)-(A2) are covered by (A3)-(A4). However, for the sake of clarity and to make the paper easy to read, we separately treat these two situations.

 Let us now present our method and their advantages when $\theta_t = 0$:  In order to establish the existence of solutions, we first only assume (A1) then we prove that  the solvability of $eq(\xi, H)$ is reduced to  the positive solvability (i.e. $Y\geq 0$) of the two BSDEs
	$eq\big(u_f(\xi_{\alpha, \beta}^+),\,  0\big)$ and $eq\big(u_f(\xi_{\alpha, \beta}^-),\,  0\big)$ (see the proof of Proposition \ref{get-gdef}). We finally show that the latter two BSDEs have simultaneously positive solutions if and only if assumption (A2) is satisfied. This shows how assumption (A2) is not a priori imposed here  but is generated along the proof.   Our method makes it possible to control more precisely the integrability condition we should impose to the  terminal value.

Here are some other advantages of the domination method: instead of $eq(\xi, H)$, we only work with the
dominating equations  $eq (\xi_1, H_1)$ and $eq (\xi_2, H_2)$ which are  more simple than the initial one $(\xi, H)$. In contrast to the papers \cite {BEK, BH1, LSM2, Kob, Tevz}, our method  allows to get the existence of  \textit{a maximal and a minimal  solution}. It moreover  allows to deal with all one dimensional BSDEs up to quadratic ones and seems unify  their treatment. Note also that, 
   the three involved terminal data $\xi$, $\xi_1$ and $\xi_2$ are not necessary integrable,
 the two dominating coefficients $H_1$, $H_2$ are merely measurable and can have arbitrarily growth. Only $H(t,y,z)$ should be continuous on $(y,\, z)$ and of  at most  quadratic growth in $z$. In return, the solutions lie in $\mathcal{C}\times\mathcal{L}^2$ and hence not necessary integrable.

 We summarize the results of some previous works in the light of assumptions (A1)-(A2). The case where $f$ is constant and $\alpha, \beta$ do not depend on $\omega$ has been considered in  \cite{Kob} where the existence of bounded solutions is established provided that the terminal value is bounded. In \cite{BH1}, the existence of solutions is obtained  when $\alpha, \beta$ and $f$ are constant ($f(y) = \frac\gamma 2$) provided that  $\exp(\gamma|\xi|e^{\beta T})$ is integrable. The authors of  \cite{BEK} consider the case where  $f$ is constant($f(y) = \frac\gamma 2$) and established the existence of solutions when $exp\big[\gamma(|\xi|+ \int_0^T\alpha_s ds)e^{\int_0^T\beta_s ds}\big]$ is integrable.  In \cite{EH2011}, the authors consider a generalized QBSDE and the function $f$ is replaced by a continuous process $r_t$,  the existence of solution is then obtained provided, when $\alpha = \beta=0 $,  that  \  $\frac{\exp(C_T|\xi|) -1}{C_T}\1_{\{C_T>0\}} + |\xi|\1_{\{C_T=0\}}$, with $C_T = \displaystyle\sup_{s\leq T}r_s$.
  The authors of \cite{beo2013, beo2017} consider the case where  $\alpha, \ \beta$ are constant and $f$ is globally integrable  and they established the  existence of solutions in $\mathcal{S}^{2}\times\mathcal{M}^{2}$ when   the terminal value is merely square integrable.

Let us briefly describe the principal methods  used in some previous papers.  When the function $f$ is constant, two methods  have been essentially developed in order to establish the existence of solutions. The first one is the monotone stability \cite{BEK, BH1, Kob, LSM2}. The second approach is based on a fixed point argument and has been introduced in  \cite {Tevz}. In the latter, the uniqueness is also obtained but it requires that the generator satisfies the so-called    Lipschitz-quadratic condition. These two methods use some a priori estimates and approximations which are sometimes difficult to obtain. It should be noted that, the papers \cite{BEK, BH1,  Kob, LSM2, Tevz} consider   the cases where  the terminal value is  bounded or at least with some exponential moments. An alternative method was recently developed in \cite{beo2013, beo2017}. This method are based on the work \cite{EH2013} where the existence of reflected QBSDEs is established without any integrability of the terminal value.  The idea,  used in \cite{beo2013, beo2017}, consists then in deriving  the existence of BSDEs from the existence of a suitable reflected BSDEs when the solutions belong to $\mathcal{S}^2\times\mathcal{M}^2$.

 We now compare our method with those of \cite{BEK, BH1, fd, Kob, LSM2, Tevz}. In the latter,  the authors proceed as follows: they first impose some integrability (or boudedness) condition on the terminal value $\xi$. Next, they establish some a priori estimates for the solutions by  using the integrability (or the boundedness)  of $\xi$.
 This allows them to prove the existence of solutions by using a suitable approximation.

	Our approach is completely different:	in order to prove the existence of solutions, we only use Lemma \ref{d} and some change of variables formulas. We do not need to establish any a priori estimates of the solutions. We do not need to construct any approximation. In contrast to 
 the previous papers, the integrability of the terminal value is not a priori imposed but obtained by solving an inverse problem. In contrast to the works \cite{BEK, BH1, fd, Kob, LSM2, Tevz},  our result covers the BSDEs with at most linear growth, and in particular it extends the results of \cite{LSM1} and \cite{ht} by taking $f=0$.


The paper is organized as follows.   In section 2, we establish the existence and uniqueness for the BSDE $(\xi, f(y)|z|^2)$ when $f$ is \textit{locally integrable} on $\R$, we also give some examples of BSDEs which have solutions without any integrability of the terminal value $\xi$.   In section 3, we begin by introducing the domination argument then we use it to establish the existence of solutions  to $eq(\xi, H)$ under conditions (A1)-(A2) and also under assumption  (A3)-(A4). Some integrability properties are also established for the solutions of $eq(\xi, H)$ under additional assumptions which will be specified below in section 3. In section 4, we treat the BSDEs with at most logarithmic growth $y\ln|y|$. Using the domination argument and some  change of variables, we show that these equations can be solved by using the quadratic BSDEs and vice-versa. Insection 5, some auxiliary results are given.

Since our approach consists in reducing the solvability  of $eq(\xi, H)$ under assumptions (A1) [resp. (A3)] to the positive solvability ($Y\geq 0$) of  $eq(u_f(\xi_{\alpha, \beta}^+), 0)$ [resp. $eq(u_f(\xi_{\alpha, \beta}^+), \theta_t|z|)$], the following two propositions which study the existence of positive solutions to these two simple BSDEs are then useful.   
 

\subsection{Two basic BSDEs} 

 The following proposition is useful in studying the solvability of $(eq(\xi, H)$ when assumption (A1) is satisfied. It characterizes the existence of positive solutions to a BSDE driven by a null  generator.   

\begin{proposition}\label{rk0} \textbf{The BSDE  $(\zeta, 0)$}.
	
	I)  According to Dudley's theorem \cite{Dud}, the following BSDE has a solution for any  $\mathcal{F}_T$--measurable random variable $\zeta$.
	\begin{equation}\label{0}
	y_t := \zeta - \int_t^T z_sdW_s
	\end{equation}
	Furthermore,

	\vskip 0.2cm $(i)$  If $\zeta$ is integrable, then equation \eqref{0} has a unique solution $(y,z)$ such that $y$ belongs to class $(D)$ given by $Y_t = \E(\zeta / \mathcal{F}_t)$. Furthermore, $z$ belongs to $\mathcal{M}^{p}$ for each $0< p < 1$.	
	
	\vskip 0.2cm
	$(ii)$ If $\zeta$ is positive and  \eqref{0} has a positive solution, then $\zeta$ is necessary integrable.
	
	\vskip 0.2cm $(iii)$ If $\zeta \neq 0$, then for any process $z$, $(0, z_t)$ could not be a solution to the BSDE $(\zeta, 0)$.

\end{proposition}

\bop \ Assertion $(i)$ can be proved by using a  usual localization and Fatou's lemma.   We prove Assertion $(ii)$.  Dudley's representation theorem allows us to show that $eq(\zeta, 0)$ has at least a solution  $(y,z)$  in $\mathcal{C}\times\mathcal{L}^2$.
Since $\zeta$ is integrable, then $y_t := \E(\zeta/\mathcal{F}_t)$ is a solution which belongs to class $(D)$. It follows that the stochastic integral $\int_{0}^{.}z_{s}dW_{s}$ is a  uniformly integrable martingale. Using Proposition 4.7, Chap. IV of \cite{RY}  (see also \cite{delyon})  and the Burkholder-Davis-Gundy inequality we show that $z$ belongs to $\mathbb{M}^{p}$, for each $0 < p < 1$. We shall prove that the process $(y,z)$ we just constructed is  actually the unique solution such that $y$ belongs to class   $(D)$. Let
$(y^{1},z^{1})$  and $(y^{2},z^{2})$ be two solutions such that $y^1$ and $y^2$ belong to class   $(D)$. It follows that $y^{1}- y^{2}$ belongs to class $(D)$ and hence the stochastic integral $(\int_{0}^{t}(z^{1}_{s} - z^{2}_{s})dW_{s})_{0\leq t \leq T}$ is a martingale in class $(D)$. It follows that $y^{1} = y^{2}$. Using the Burkholder-Davis-Gundy inequality, we show that $\E\left[\left(\int_{0}^{T}|z^{1}_{s} - z^{2}_{s}|^2\right)^{p/2}\right] = 0$.   Assertion $(ii)$ is proved. We shall prove $(iii)$. Let $\zeta \neq 0$. Assume that there exists a process $z$ such that $(0,z)$ is a solution to $eq(\zeta, 0)$. Then for any $t \leq T$, \ $0 = \zeta - \int_t^T  z_s dW_s$ which implies that $\zeta =0$ by putting $t= T$. 
 \eop


\vskip 0.2cm The following proposition, which is taken from \cite{EH2011} (Proposition 6.1 of \cite{EH2011}), gives a necessary and sufficient condition which ensures the existence of  positive solutions to the  BSDEs driven by the generator $ \theta_t|z|$. This proposition is useful when assumption (A3) is satisfied.

\begin{proposition} \label{propehz} ( \cite{EH2011}, Proposition 6.1). \textbf{The BSDE $(\zeta, \ \theta_t|z|^2)$}.     Let $\zeta$ be a positive  $\mathcal{F}_T$--measurable random variable. The
	BSDE
	\begin{equation} \label{ehz}
	y_{t}={\zeta} +\int_t^T \theta_s |z_s| ds
	-\int_t^T z_{s}dW_{s}\,,  \qquad t\leq T
	\end{equation}
	has a positive solution  if and only if
	\begin{equation}\label{ehxi}
	\sup_{\pi\in\sum} \E\left(\Gamma_{0, T}^{\pi}\zeta\right) :=\sup_{\pi\in\sum} \E\left(e^{\int_0^T \theta_u \pi_udW_u -\frac12 \int_0^T
		\theta_u^2 |\pi_u|^2du}{\zeta}\right)  < +\infty
	\end{equation}	
	where   \
	$ \sum := \left\{  \pi\in {\cal L}^{2},\, |\pi|\in  \{0, 1\},  \  a.e. \ {\rm{and}} \  \displaystyle \  ess\sup_{\omega} \int_0^T
	\theta_u^2 |\pi_u|^2du < +\infty\right\}  $	
	
	\vskip 0.2cm
	In this case, there exist $\bar{z}\in {\cal L}^{2}$ and $\bar{y}_t := \displaystyle{ess\sup_{\pi\in \sum}} \E(\Gamma_{t,
		T}^{\pi}{\zeta}|{\cal F}_t)
	$
	such that $(\bar{y},
	\bar{z})$
	is the minimal solution of Equation (\ref{ehz}). 	
	Furthermore, $\bar{y}_t \geq \E({\zeta}|{\cal
		F}_t)\geq 0$,\,\ for each $t\in [0,T]$.

	
\end{proposition}





	

\section{The BSDE$(\xi, f(y)|z|^2)$}

 The BSDE$(\xi, f(y)|z|^2)$ will be used in order to solve the general equation $(\xi, H)$ with $|H(t, y, z)| \leq \alpha_t|y|+ \beta_t|y| +\theta_t|z| + f(|y|)|z|^2$. However, since $eq(\xi, f(y)|z|^2)$ is interesting itself and do not need the domination argument, 	
 we give in this subsection a complete study of this equation in the case where the function $f$ is locally integrable but not necessary continuous. A characterization of the existence of solution is given for this equations. This is related to the function $u_f$, and for instance, when $u_f(\R) = \R$ then the BSDE$(\xi, f(y)|z|^2)$ has a solution for each terminal value $\xi$. No integrability is required to the terminal value $\xi$. We start this section by some examples which are covered by the present work and not covered by those of \cite{beo2013, beo2017, BEK, BH1, EH2011, fd, Kob, LSM2, Tevz}.


\subsection{Some examples of QBSDEs with non integrable terminal value}

	\begin{example}\label{ex1}
		Let $f_1$ be a bounded function which is globally integrable on $\mathbb{R}$. We assume that $f$ is bounded by 1 for simplicity. Clearly, the generator $H_1(y, z) := f_1(y)|z|^2$ satisfies $|H(y, z)| \leq |z|^2$. Hence, $H_1(y, z)$ is of quadratic growth.  It was shown in \cite{beo2013, beo2017} that the QBSDE $(\xi, f_1(y)|z|^2)$  has a solution without any integrability condition on the terminal value $\xi$. Moreover, when $\xi$ is square integrable then $eq(\xi, f_1(y)|z|^2)$ has a unique solution in $\mathcal{S}^{2}\times \mathcal{M}^{2}$.
	\end{example}
	Indeed, for a given locally integrable function $f$, the transformation $u_f$ defined in \eqref{u} is a one to one function from $\mathbb{R}$ onto $\mathbb{R}$. Both $u_{f}$ and its inverse belong to the Sobolev space ${W}_{1,\, loc}^{2}$. Using Itô-Krylov's formula for BSDEs (see \cite{beo2013, beo2017}), it follows that
	$eq(\xi ,f_1(y)|z|^2)$ has a solution  if and only if $eq(u_{f_1}(\xi ),0)$ has a solution.
	Since $f_1$ is globally integrable, then $u_{f_1}$ and its inverse are uniformly Lipschitz. It follows that $u_{f_1}(\xi)$ is square integrable if and only if $\xi$ is square integrable.
	Therefore,   $eq(\xi ,f_1(y)|z|^2)$ has a unique solution in $\mathcal{S}
	^{2}\times \mathcal{M}^{2}$ whenever $\xi$ is merely square integrable.
	This also shows that the convexity of the generator is not necessary to the uniqueness.

	\begin{example}\label{exex4} The functions \
		\begin{equation}\label{f34}
		f_2(y) := e^y \qquad\qquad and  \qquad\qquad f_3(y) := \1_{\{y > 0 \}} + \1_{\{y \leq 0 \}}e^y
		\end{equation}
		are not globally integrable on $\mathbb{R}$.
		However, the same argument shows that  $eq(\xi, f_i(y)|z|^2)$ has a solution without any integrability of $\xi$, for $i \in\{ 2, 3\}$. Note that $f_2$ is neither   globally integrable nor bounded.
	\end{example}



	\subsection{ Existence  of BSDE$(\xi, f(y)|z|^2)$ with $f$ locally integrable}

	 Let $f:\R \longrightarrow \R$ be a locally integrable function. The goal of this part is to explain how some $eq(\xi, f(y)|z|^2)$ has solutions without integrability of $\xi$ and others need the integrability of $\xi$.
		
		Consider the QBSDE
		\begin{equation}\label{fz2}
		Y_{t}=\xi +\int_{t}^{T}f(Y_s) |Z_{s}|^2ds-\int_{t}^{T}Z_{s}dW_{s},\, \ \ 0\leq
		t\leq T.
		\end{equation}
		
		The function $u_{f}$ defined by \eqref{u} belongs to the Sobolev space $W^2_{1, loc}(\R)$. Therefore, applying Itô-Krylov's formula to $u_{f}$ (see \cite{beo2013, beo2017}), one can show that $(Y,Z)$ is a solution to equation \eqref{fz2}  if and only if $(\bar{Y},\bar{Z}):= (u_{f}(Y), u_{f}'(Y)Z) $ is a solution to the BSDE
		
		\begin{equation}\label{zero}
		\bar{Y}_{t}=u_{f}(\xi) -\int_{t}^{T}\bar{Z}_{s}dW_{s},\, \ \ 0\leq
		t\leq T .
		\end{equation}
		
		According to Dudley's representation theorem, $eq(u_f(\xi), 0)$ has a solution for any $\mathcal{F}_T$-measurable random variable $u_f(\xi) $. No integrability is required to  $u_f(\xi)$. But, our problem is to solve $eq(\xi, f(y)|z|^2)$. This is the subject of the following proposition.

		\begin{proposition}\label{quadraticpure}
		Let $f$ be a locally integrable function  and $\xi $ a  $\mathcal{F}_T$-measurable random variable.

		 $(i)$ \ If $u_f(\R) = \R$, then   $eq(\xi, f(y)|z|^2)$ has a solution. No integrability is needed for $\xi$.
			
			 $(ii)$ \ Let $u_f(\R) \neq \R$. If $f$ is positive, $u_f(y)$ is then increasing and  $\displaystyle\lim_{y\rightarrow\infty}u_f(y) = + \infty$.
			\ Assume that $\displaystyle\lim_{y\rightarrow -\infty}u_f(y) =  c > - \infty$. Then, necessary $\bar u_f(\xi) := u_f(\xi) - c$ is  integrable. In this case, $eq(\xi, f(y)|z|^2)$ has a solution given by $Y_s = \bar u_f^{-1}(\E[\bar u_f(\xi)/\mathcal{F}_s])$.  The case   $f$ negative goes similarly.

		\end{proposition}
			
			\bop Assertion $(i)$ is simple. We prove assertion $(ii)$.
			 The function  $\bar u_f(y) := u_f(y) - c$ \ belongs to the Sobolev space  ${W}^2_{1, loc}(\mathbb{R})$ and it is one to one from $\mathbb{R}$ into $\mathbb{R}_+$.  Hence, Itô-Krylov's formula applied to $\bar u_{f}$ shows that $(Y,Z)$ is a solution to equation \eqref{fz2}  if and only if $(\bar{Y},\bar{Z}):= (\bar u_{f}(Y),  \bar u_{f}'(Y)Z) $ is a solution to the BSDE
			\begin{equation}\label{vzero}
			\bar{Y}_{t}=\bar{u}_{f}(\xi) -\int_{t}^{T}\bar{Z}_{s}dW_{s},\, \ \ 0\leq
			t\leq T.
			\end{equation}
			Since $\bar{Y}:= \bar u_{f}(Y) $ is positive,  Therefore, $eq(\xi, f(y)|z|^2)$ has a solution if and only if $\bar u_f(\xi)$ is integrable. We then deduce that for any $s\leq T$, \ $Y_s = \bar u_f^{-1}(\E[\bar u_f(\xi)/\mathcal{F}_s])$.


	\subsection{Uniqueness and comparison for $eq(\xi, f(y)|z|^2)$ with $f$ locally integrable}
\begin{proposition}\label{uniqfy}
Let $f:\R \longrightarrow \R$ be locally integrable. Let $u_f$ be the function defined in Lemma \ref{zvonkinM}-I, $v$ the function defined in Lemma \ref{zvonkinM}-II and  $w$ the fonction   defined in Lemma \ref{zvonkinM}-III. Assume that $u_f(\xi)$ is integrable. Then, the BSDE$(\xi, f(y)|z|^2)$ has a unique solution $(Y,Z)$ such that $u_f(Y)$ belongs to class $(D)$.

\noindent If moreover,

(i) \  $v(Y)$ belongs to class (D), then  \  $\E\int_0^T|Z_s|^2ds < \infty$.

(ii) \ $w(Y)$ belongs to class (D), then  \  $\E\int_0^T|f(Y_s)||Z_s|^2ds < \infty$.
\end{proposition}

\bop 	\
 The BSDE$(\xi, f(y)|z|^2)$  has a solution if and only if the BSDE $(u_f(\xi), 0)$ has a solution. But $eq(u_f(\xi), 0)$  has a solution by Dudley's representation theorem. This gives the existence of solutions. We prove the uniqueness.  Let $Y^1$ and $Y^2$ be two solutions of $eq(\xi, f(y)|z|^2)$ such that $u_f(Y^1)$ and $u_f(Y^2)$ belong to class $(D)$.  Arguing as in the proof of Proposition  \ref{rk0}, we show that $u_f(Y^1) = u_f(Y^2)$ which implies that $Y^1= Y^2$ since $u_f$ is one to one.  Arguing again as in the proof of Proposition \ref{rk0}, we show that $\int_0^T|Z^1_s - Z^2_s|^2ds = 0$ \ $a.s$, since $u_f'$ is strictly positive.

We prove $(i)$.
Let $v$ be the map defined in Lemma \ref{zvonkinM}-II).
For $N>0$, let
$\tau _{N}:=\inf \{t>0:|Y_{t}|+\int_{0}^{t}|v^{\prime
}(Y_{s})|^{2}|Z_{s}|^{2}ds\geq N\}\wedge T$. Set $\mbox{sgn}(x)=1$ if $x\geq
0 $ and $\mbox{sgn}(x)=-1$ if $x<0$. Since the map $x\mapsto v(|x|)$ belongs to $W^{2}_{1,\, loc}(
\mathbb{R})$, then thanks to  Itô-Krylov's formula for BSDEs (see \cite{beo2013, beo2017}), we have for any $t\in [0,T]$,
\begin{align*}
v(|Y_{0}|) & = v(|Y_{t\wedge \tau _{N}}|)
+ \int_{0}^{t\wedge \tau _{N}}\!\!\left[
\mbox{sgn}(Y_{s})v^{\prime }(|Y_{s}|)f(Y_{s})|Z_{s}|^{2}-\frac{1}{2}v^{\prime
	\prime }(|Y_{s}|)|Z_{s}|^{2}\right]\!\!ds
\\
& \quad \ - \int_{0}^{t\wedge \tau _{N}}\mbox{sgn}(Y_{s})v^{\prime
}(|Y_{s}|)Z_{s}dW_{s} \, .
\end{align*}
Lemma \ref{zvonkinM}-II)  allows us to  deduce that for any $N>0$,

\begin{align}\label{estimz}
\frac{1}{2}\mathbb{E}\!\int_{0}^{t\wedge \tau _{N}}\!\!|Z_{s}|^{2}ds
&\leq
\E\big[v(|Y_{t\wedge \tau _{N}}|)\big]
\\
&\leq
\sup_{\tau\leq T}\E\big[v(|Y_{ \tau }|)\big]
\end{align}
where the supremum in the first right hand term, is taken over all stopping times $\tau\leq T$.

Since the process $\left[v(|Y_{t}|)\right]$ belongs to class $(D)$, the proof is completed by  using Fatou's lemma.

We  prove $(ii)$.   Without loss, we assume that $f$ is positive. Let $w$ be the map defined in Lemma \ref{zvonkinM}-III).
For $N>0$, let
$\tau _{N}:=\inf \{t>0:|Y_{t}|+\int_{0}^{t}|w^{\prime
}(Y_{s})|^{2}|Z_{s}|^{2}ds\geq N\}\wedge T$. Set $\mbox{sgn}(x)=1$ if $x\geq
0 $ and $\mbox{sgn}(x)=-1$ if $x<0$. Since the map $x\mapsto w(|x|)$ belongs to $W^{2}_{1,\, loc}(
\mathbb{R})$, then thanks to  Itô-Krylov's formula for BSDEs (see \cite{beo2013, beo2017}), we have for any $t\in [0,T]$,
\begin{align*}
w(|Y_{0}|) & = w(|Y_{t\wedge \tau _{N}}|)
+ \int_{0}^{t\wedge \tau _{N}}\!\!\left[
\mbox{sgn}(Y_{s})w^{\prime }(|Y_{s}|)f(Y_{s})|Z_{s}|^{2}-\frac{1}{2}w^{\prime
	\prime }(|Y_{s}|)|Z_{s}|^{2}\right]\!\!ds
\\
& \quad \ - \int_{0}^{t\wedge \tau _{N}}\mbox{sgn}(Y_{s})w^{\prime
}(|Y_{s}|)Z_{s}dW_{s} .
\end{align*}
Assumption \textbf{(A1)} and
Lemma \ref{zvonkinM}-III)  allow us to  show that for any $N>0$,

\begin{align}\label{estimz}
\frac{1}{2}\mathbb{E}\!\int_{0}^{t\wedge \tau _{N}}\!\!f(Y_{s})|Z_{s}|^{2}ds
&\leq
\E[w(|Y_{t\wedge \tau _{N}}|)]
\end{align}
where the supremum in the first right hand term, is taken over all stopping times $\tau\leq T$.

\noindent Since $w(Y)$ belongs to class $(D)$, the proof is completed by sending $N$ to $\infty$ and using Fatou's lemma. Proposition  \ref{uniqfy} is proved. \eop


\begin{remark}\label{integrbfyzcarre}
$(i)$ \ Propositions \ref{quadraticpure} and \ref{uniqfy} give the existence of solution in $\mathcal{C}\times\mathcal{L}^2$ with and without integrability of $u_f(\xi)$. If one wants to get  more integrability of the solutions, then it is  sufficient  to impose more integrability on $u_f(\xi)$. For instance, when $\xi$ is bounded then the solution $(Y,Z)$ is such that $Y$ is bounded and $(\int_0^t Z_sdW_s)_{0\leq t \leq T}$ is a BMO martingale.

$(ii)$ \ Note that, in contrast to \cite{beo2013, beo2017}, the integrability of $u_f(\xi)$ is not equivalent to that of  $\xi$.
\end{remark}


	
\begin{proposition}
	\label{ct}{\rm (Comparison)} Let $\xi_1$, $\xi_2$ be $\mathcal{F}_T$--measurable. Let $f_1$, $f_2$ be elements of \   $\mathbb{L}_{loc}^{1}(
	\mathbb{R})$. Assume that $u_{f_1}(_xi_1)$ and $u_{f_2}(_xi_2)$ are integrable.  Let $( Y^{f_1},Z^{f_1})$, $( Y^{f_2},Z^{f_2}) $ be respectively the unique 
	solution in class $(D)$ of $eq(\xi_1, f_1(y)|z|^2)$ and $eq(\xi_2, f_2(y)|z|^2)$.
	Assume that $\xi_{1}\leq\xi_{2}$ a.s. and $f_1 \leq f_2$ $a.e$. Then $
	Y_{t}^{f_1}\leq Y_{t}^{f_2}$ for all $t$ $\mathbb{P}$--a.s.
\end{proposition} 


\bop According to Proposition \ref{uniqfy}, the solutions $ Y^{f_1}
$ and $Y^{f_2} $ belong to class $(D)$. Arguing as in the proof of Proposition \ref{uniqfy}, one can show that the processes $\big(\int_0^t u'_{f_1}(Y_{s}^{f_1})Z_{s}^{f_1}dW_s\big)_{0\leq t\leq T}$ and $\big(\int_0^t u'_{f_2}(Y_{s}^{f_2})Z_{s}^{f_2}dW_s\big)_{0\leq t\leq T}$ are uniformly integrable martingales. Using the Burkhölder-Davis-Gundy inequality and the fact that $f_1\leq f_2$, we show that the process $\big(\int_0^t u'_{f_1}(Y_{s}^{f_2})Z_{s}^{f_2}dW_s\big)_{0\leq t\leq T}$ is a uniformly integrable martingale. The rest of the proof can be  performed as that of Proposition 3.2. in   \cite{beo2017}.  \eop

	\section{ BSDE$(\xi, H)$ }

To deal with more general BSDEs, we need the domination argument which will be present in the following subsection.

\subsection{The domination argument}

The domination argument implicitly appears in \cite{beo2013, beo2017} in a particular situation where the two dominating solutions belong to $\mathcal{S}^2\times\mathcal{M}^2$, the function $f$ is globally integrable on $\R$, the three terminal values are square integrable, the two dominating coefficients $H_1$, $H_2$ are of  quadratic growth in $z$ and of linear growth in $y$. Here, this argument is presented in a general framework, that is:   the two dominating  solutions lie in $\mathcal{C}\times\mathcal{L}^2$, the three involved terminal data $\xi$, $\xi_1$ and $\xi_2$ are not necessary integrable,
the two dominating coefficients $H_1$, $H_2$ are merely measurable and can have arbitrarily growth. Only $H(t,y,z)$ should be continuous on $(y,\, z)$ and of  at most  quadratic growth in $z$.

\begin{definition}\label{defdom} (Domination conditions) We say that  the data $(\xi, H)$  satisfy a domination condition if  there exist two $(\mathcal{F}_t)$ progressively measurable processes $H_1$ and  $H_2$, two $(\mathcal{F}_T)$ measurable random variables $\xi_1$ and $\xi_2$ such that:
	\vskip 0.1cm\noindent
	\rm{\textbf{(D1)}}  \quad $\xi_1 \leq  \xi \leq \xi_2$
	\vskip 0.1cm\noindent
	\rm{\textbf{(D2)}} \quad   $eq(\xi_1, H_1)$ and $eq(\xi_2, H_2)$ have two solutions $(Y^1,Z^1)$ and $(Y^2,Z^2)$ such that:
	\vskip 0.1cm\noindent
	\hskip 0.7cm\rm{\textbf{(a)}}     \quad $Y^1 \leq Y^2$,
	\vskip 0.1cm\noindent
	\hskip 0.7cm\rm{\textbf{(b)}} \quad for every $
	(t,\omega)$,  $y\in [Y_t^1(\omega), \ Y_t^2(\omega)]$ and $z\in \mathbb{
		R}^d$,
	\vskip 0.1cm\noindent
	\hskip 1,5cm $(i)$ \ $H_1(t,y,z) \leq  H(t,y,z) \leq H_2(t,y,z)$
	\vskip 0.1cm\noindent
	\hskip 1,5cm $(ii)$ \
	$
	|H(t, \omega, y, z )| \leq \eta_t(\omega)+ C_t(\omega)|z|^2
	$
	\vskip 0.1cm\noindent
	\textsl{where  $C$ and $\eta$ are $\mathcal{F}_t$-adapted processes such that $C$ is continuous and $\eta$ satisfies for each $\omega$, \  $\int_0^T|\eta_s(\omega)|ds < \infty$.}
\end{definition}

\begin{lemma}\label{d}(Existence by domination)
	Let $H$ be continuous in $(y,z)$ for $a.e. \ (t,\omega)$. Assume moreover that $(\xi, H)$  satisfy the domination conditions (D1)--(D2). Then,
	
	$(i)$ The BSDE $(\xi, H)$ has at least one solution $(Y,Z)$ such that
	$Y^1 \leq Y \leq Y^2$.
	
	$(ii)$ Among all solutions which lie  between 	$Y^1$ and  $ Y^2$, there exist a maximal and a minimal solution.
\end{lemma}

This lemma  directly gives the existence of solutions. We do not need any a priori estimates nor approximation.   The idea of the proof consists in deriving the existence of solutions for the BSDE
without reflection from solutions of a suitable QBSDE with two reflecting barriers. To this end, we use the remarkable result of Essaky $\&$ Hassani
(\cite{EH2013}, Theorem 3.2) which establishes the existence of solutions for
reflected QBSDEs without assuming any integrability condition on the
terminal value. For the self-contained, this result is stated in Theorem \ref{EH2013} in Appendix. 


\vskip 0.3cm\noindent
\textbf{Proof of Lemma \ref{d}}  Using Theorem 3.2 in \cite{EH2013} (see Theorem \ref{EH2013} in Apendix) with $L
=Y^{H_1}$\ and \ $ U=Y^{H_2}$,  there exists a process $(Y,Z, K^+, K^-)$  such that $(Y,Z)$ belongs to $\mathcal{C}\times\mathcal{L}^2$ and $(Y,Z, K^+, K^-)$ satisfies the following reflected BSDE, for $t\in [0, \ T]$,

\begin{equation}  \label{RBSDEHdominé}
\left\{
\begin{array}{ll}
(i) & Y_{t}=\xi +\displaystyle\int_{t}^{T}H(s,Y_{s},Z_{s})ds
-\displaystyle
\int_{t}^{T}Z_{s}dB_{s}\, \\
& \qquad \quad +\displaystyle
\int_{t}^{T}dK_{s}^+ -\displaystyle\int_{t}^{T}dK_{s}^- \\
(ii) & \forall \ t\leq T,\, \quad Y_t^{H_1} \leq Y_{t}\leq Y_t^{H_2}, \\
(iii) & \displaystyle\int_{0}^{T}( Y_{t}-Y_t^{H_1}) dK_{t}^+ =
\displaystyle
\int_{0}^{T}( Y_t^{H_2}-Y_{t}) dK_{t}^-=0,\,\, \mbox{a.s.}, \\
(iv) & K_0^+ =K_0^- =0, \,\,\,\, K^+, K^- \,\,
\mbox{are continuous
	nondecreasing}, \\
(v) & dK^+ \bot dK^- .
\end{array}
\right.
\end{equation}
Moreover, equation \eqref{RBSDEHdominé} has a minimal solution and a maximal solution.

It remains to show that $dK^+ = dK^- = 0$.  Since $Y_t^{H_2}$ is a
solution to the BSDE $ eq(\xi_2, H_2)$, then Tanaka's formula
applied to $(Y_t^{H_2}-Y_t)^+$ shows that
\begin{align*}
(Y_{t}^{H_2}-Y_{t})^{+} & =(Y_{0}^{H_2}-Y_{0})^{+}+\int_{0}^{t} \1
_{\{Y_{s}^{H_2} > Y_s\}} [H(s,Y_s,Z_s)-H_2(s,Y_s^{H_2},Z_s^{H_2})]ds \\
& \ + \int_{0}^{t} \1_{\{Y_{s}^{H_2} > Y_s\}}(dK_{s}^+-dK_{s}^-) +
\int_{0}^{t} \1_{\{Y_{s}^{H_2} > Y_s\}}(Z_{s}^{H_2}-Z_s)dW_s \\
& \ + L_t^0(Y^{H_2}-Y)
\end{align*}
where $L_t^0(Y^{H_2}-Y)$ is the local time at time $t$ and level $0$ of the
semimartingale $(Y^{H_2}-Y)$.

\vskip 0.15cm\noindent Since  $(Y_{t}^{H_2}-Y_{t})^{+} =
(Y_{t}^{H_2}-Y_{t})$,  then by identifying the terms of $
(Y_{t}^{H_2}-Y_{t})^{+}$ with those of $(Y_{t}^{H_2}-Y_{t})$ one can show
that:
\begin{align*}
(Z_s - Z_s^{H_2})\1_{\{Y_{s}^{H_2}=Y_{s}\}} = 0 \quad \text{for} \ a.e. \ (s, \omega)
\end{align*}
and
\begin{align*}
\int_{0}^{t}\1_{\{Y_{s}^{H_2}=Y_{s}
	\}}(dK_{s}^{+}-dK_{s}^{-})=L_t^0(Y^{H_2}-Y) + \int_{0}^{t}\1
_{\{Y_{s}^{H_2}=Y_{s} \}}[H_2(s,Y_{s}^{H_2},Z_{s}^{H_2})-H(s,Y_{s},Z_{s})]ds
\end{align*}
Since $\displaystyle\int_{0}^{t}\1_{\{Y_{s}^{H_2}=Y_{s}\}}dK_{s}^{+}=0$, it holds
that
\begin{align*}
0\leq & L_t^0(Y^{H_2}-Y)+\int_{0}^{t}\1_{\{Y_{s}^{H_2}=Y_{s}
	\}}[H_2(s,Y_{s}^{H_2},Z_{s}^{H_2})-H(s,Y_{s},Z_{s})]ds = -\int_{0}^{t}\1
_{\{Y_{s}^{H_2}=Y_{s}\}}dK_{s}^{-}\leq 0
\end{align*}
Hence,
$\int_{0}^{t}\1_{\{Y_{s}^{H_2}=Y_{s}\}}dK_{s}^{-}=0$, which
implies that \ $dK^{-}=0$. Similar arguments allow to show that
$ dK^{+}=0 $. Therefore $(Y,Z)$ is a solution to the (non reflected)
BSDE $ eq(\xi ,H)$.
\eop

As byproducts of the domination lemma, we  establish in this section the existence of solutions to $eq(\xi, H)$ first when (A1)-(A2) are satisfied and next  when (A3)-(A4)  hold. Some integrability properties of the solutions are also established under additional assumptions which will be specified below.

	\subsection{ BSDE$(\xi, H)$ with $|H(t, y, z)|   \leq \alpha_t + \beta_t\vert y\vert  + f(|y|)\vert z\vert^{2}$}

The goal is to solve $eq(\xi, H)$. If one try to follow the proofs given in \cite{BEK, BH1}, we should establish some a priori estimates of solutions then use a suitable approximation. This way is not efficient in our situation and in particular the exponential transformation can not be applied because $f$ is not constant. Furthermore, in contrast to the domination argument, the method of \cite{BEK, BH1} can not be applied when the integrability of the terminal value is not a priori fixed.   For the same reason, the argument developed in \cite{fd, DHR1, Tevz} are also not effective in our situation.  Note moreover that the  methods used in \cite{BEK, BH1, fd, DHR1, Tevz} can not allow to prove the existence of a maximal and a minimal  solution. The method used in \cite{beo2013, beo2017} does not work in our situation since $f$ is not globally integrable. The questions which then arise are : what condition we should impose to the terminal value $\xi$ in order to get the existence of solution when $f$ is not globally integrable? How do we proceed in this case ?  This is the subject of the next subsection.
 To prove the existence of solutions to $eq(\xi, H)$, our strategy consists in using Lemma \ref{d} which allows to work without any a priori  integrability condition on the terminal value $\xi$.  Therefore, we start by assuming that only (A1) is satisfied.
Let
\begin{equation}\label{g}
g(t, y, z) := \alpha_t + \beta_t\vert y\vert  + f(|y|)\vert z\vert^{2}
\end{equation}
According to Lemma \ref{d}, to establish the existence of solutions to $eq(\xi, H)$, it is enough to show that $eq(\xi^+, g)$ has a solution $(Y^g, Z^g)$ and $eq(-\xi^-, -g)$ has a solution $(Y^{-g}, Z^{-g})$ such that
\begin{equation}\label{y(-g)yg}
Y^{-g}\leq Y^g.
\end{equation}
In  \cite{beo2013, beo2017}, the fact that $\xi$ is assumed square integrable and $f$ is globally integrable make simple the solvability of $eq(\xi^+, g)$ and $eq(-\xi^-, -g)$ in $\mathcal{S}^2\times\mathcal{M}^2$ from which we easly deduce inequality \eqref{y(-g)yg}.
Question:  

\textit{how to prove  inequality \eqref{y(-g)yg} when we  do not have any  information  on the integrability of $\xi$ nor on the integrability of the solutions? }

We emphasize that the comparison theorem does not work in this situation. But we need to prove inequality \eqref{y(-g)yg} in order to establish the existence of solutions to  $eq(\xi, H)$.   We proceed as follows :	  we  assume  that (A1) holds, \textit{we force $Y^{-g}$ to be negative and $Y^{g}$ to be positive  then we deduce the integrability condition [namely (A2)] which we should impose to the terminal value $ \xi $ in order to get the existence of solution}.  Hence, assumption (A2) is generated by solving  an \textit{inverse problem}. 
Lemma \ref{d} plays a key role in our proof. In particular, it allows us to reduce the solvability of $eq(\xi^+, g)$ to that of $eq(\xi^+_{\alpha, \beta}, 0)$ from which we derive assumption (A2).

	\subsubsection{ The BSDE $(\xi, \alpha_t + \beta_t |y| + f(|y|)|z|^2)$}

The goal of this subsection is to show that the BSDE $(\xi^+, g)$ has a positive solution and the BSDE $(\xi^-, -g)$ has a negative solution.
The method used in \cite{beo2013, beo2017} to  show this consists in reducing the solvability of these two BSDEs to the solvability of two BSDEs with linear growth. These computations do not work in our situation, since $f$ is not globally integrable. The proof of following proposition shows how to get the existence of a  positive to $eq(\xi^+, g)$ by using the domination argument when we do not have information on the integrability terminal value $\xi$. It also allows the determine the integrability of we should impose to $\xi$ by solving an inverse problem.

	\begin{proposition}\label{get-gdef}
		Let  $\xi$ be an $\mathcal{F}_T-$measurable random variable. Let $f$ be as in assumption (A1) and $g$ be the function defined by \eqref{g}.
		
		$(i)$  \ The BSDE \ $(\xi^+, g)$ has a positive solution
		\  if \  the BSDE \  $\big(u_{f}(\xi_{\alpha, \beta}^+),\ 0\big)$ has a positive solution $(Y, Z)$ such that \
		$Y \geq u_{f}\left[e^{\int_0^T\beta_s ds}\left(\int_0^T\alpha_s ds\right)\right]$.
		
		$(ii)$  \ The BSDE \ $(-\xi^-, -g)$ has a negative solution  \ if \ the BSDE
		\\  $\big(u_{f}(\xi_{\alpha, \beta}^-),\ 0\big)$ has a positive solution $(Y, Z)$ such that \
		$Y \geq u_{f}\left[e^{\int_0^T\beta_s ds}\left(\int_0^T\alpha_s ds\right)\right]$.
		
		$(iii)$ If moreover assumption (A2) is satisfied, then $eq(\xi^+, g)$ has a positive solution and $eq(\xi^-, -g)$ has a negative solution.
	\end{proposition}
	
	
	\bop \  $(i)$ \
	For the simplicity of notations, we assume that $\alpha$ and $\beta$ are constant. Note that $(Y,Z)$ is a positive solution to $eq(\xi^+, g)$ if and only if  $(Y,Z)$ is a positive solution to the BSDE
	\begin{equation}\label{eqg1}
	Y_{t}=\xi^+ +\int_{t}^{T} \alpha + \beta Y_s + f\left(Y_{s}\right) |Z_{s}|^{2}ds-\int_{t}^{T}
	Z_{s}dW_{s}.
	\end{equation}
	Therefore, it is enough to prove that: if \  $eq(u_f[e^{\beta T}(\xi^+ + \alpha T)],\ 0)$ has a  solution $(Y^{0}, Z^{0})$ such that $Y^{0} \geq u_f( \alpha T e^{\beta T})$ then equation \eqref{eqg1} has a positive solution.

	Return back to BSDE \eqref{eqg1}. By putting  $(Y^1_t, Z^1_t):= (Y_t+ \alpha t, Z_t)$, we see that equation  \eqref{eqg1} has a positive solution if and only if the BSDE
	\begin{equation}\label{fY1}
	Y^1_{t}=\xi^+ + \alpha T +\int_{t}^{T}  \beta (Y^1_s-\alpha s) + f\left(Y^1_{s}-\alpha s)\right) |Z^1_{s}|^{2}ds-\int_{t}^{T}
	Z^1_{s}dW_{s}
	\end{equation}
	has a solution $(Y^1, Z^1)$ such that for each t, $Y^1_t \geq \alpha t$.
	
	Consider now the BSDE
	\begin{equation}\label{alphat}
	Y_{t}= \alpha T -\int_{t}^{T}Z_{s}dW_{s}
	\end{equation}
	Clearly
	\begin{itemize}

	\item $(Y, Z) = (\alpha T, 0)$ is a solution to the BSDE \eqref{alphat},
	
	\item $0 \leq \xi^+ + \alpha T \leq \xi^+ + \alpha T$,
	
	\item for any  $y\geq \alpha T$,  \ $0 \leq \beta (y-\alpha s) + f\left(y-\alpha s\right) |z|^{2} \leq \beta y + f\left(y\right) |z|^{2}$ since $f$ is increasing.

\end{itemize}
	
	\noindent
	Therefore, using Lemma \ref{d} [with $\xi_1 = 0$, $\xi= \xi^+ + \alpha T$ = $ \xi_2$,  $H_1 = 0$ and $H_2 = \beta y + f(y)|z|^2$], we show that equation \eqref{fY1} has a solution $(Y^1, Z^1)$ satisfying $Y^1\geq \alpha T$ if the BSDE
	\begin{equation}\label{fY2}
	Y^2_{t}=\xi^+ + \alpha T +\int_{t}^{T}  \beta (Y^2_s) + f\left(Y^2_s\right) |Z^2_{s}|^{2}ds-\int_{t}^{T}
	Z^2_{s}dW_{s}
	\end{equation}
	has a solution $(Y^2, Z^2)$ satisfying $Y^2\geq \alpha T$.
	
	But, Itô's formula shows that $(Y^2, Z^2)$ is a solution to equation \eqref{fY2} satisfying $Y^2\geq \alpha T$ if and only if the process $(Y^3_t, Z^3_t) := (Y^2_te^{\beta t}, Z^2_te^{\beta t})$ is a solution to the BSDE
	
	\begin{equation}\label{fY3}
	Y^3_{t}=(\xi^+ + \alpha T)e^{\beta T} +\int_{t}^{T}  f\left(Y^3_se^{-\beta s})\right) |Z^3_{s}|^{2}e^{-\beta s}ds-\int_{t}^{T}
	Z^3_{s}dW_{s}
	\end{equation}
	satisfying $Y^3 \geq \alpha T e^{\beta T}$. 
	
	Since $f$ is increasing and continuous then, as previously, we use again  Lemma \ref{d} [with $\xi_1 = 0$, $\xi= (\xi^+ + \alpha T)e^{\beta T}$ = $ \xi_2$,  $H_1 = 0$ and $H_2 = f(y)|z|^2$], to  show that equation \eqref{fY3} has a solution $(Y^3, Z^3)$ such that  $Y^3\geq \alpha T e^{\beta T}$ if  $eq([\xi^+ + \alpha T]e^{\beta T}, f\left(y\right) |z|^{2})$ has a solution $(Y^4, Z^4)$ satisfying $Y^4 \geq \alpha T e^{\beta T}$. Applying Itô's formula to $u_f(Y^4_t)$, we show that  $eq([\xi^+ + \alpha T]e^{\beta T}, f\left(y\right) |z|^{2})$ has a solution $(Y^4, Z^4)$ such that $Y^4 \geq \alpha T e^{\beta T}$ if and only if  $eq(u_f[\xi^+ + \alpha T]e^{\beta T}), 0)$ has a solution $(Y^5, Z^5)$ satisfying $Y^5 \geq u_f(\alpha T e^{\beta T})$. According to Proposition \ref{rk0}, the latter is equivalent to the fact that $u_f[e^{\beta T}(\xi^+ + \alpha T)]$ is integrable. Assertion $(i)$ is proved.
	
	Assertion $(ii)$ can be proved similarly, since  $(Y,Z)$ is a negative solution to $eq(-\xi^-, -g)$ if and only if  $(Y',Z'):=(-Y,-Z)$ is a positive solution to the BSDE
	\begin{equation}\label{eqg'}
	Y'_{t}=\xi^-\int_{t}^{T} \alpha + \beta Y_s' + f\left(Y_{s}'\right) Z_{s}^{2}ds-\int_{t}^{T}
	Z_{s}'dW_{s}
	\end{equation}
	 Proposition \ref{get-gdef} is proved. \eop

 \subsubsection{ The BSDE $(\xi, H)$ with $|H(t,y,z)| \leq \alpha + \beta |y| + f(|y|)|z|^2$}




The following theorem is deduced from  Proposition \ref{get-gdef} and Lemma \ref{d}. It  covers the previous results established in \cite{BEK, BH1, DHR1, fd, Kob} and many others situations which are not covered by the previous works on QBSDE. For instance, we cover the cases: \
$H(y,z) = y|z|$ and also $H(y,z) =  \alpha_t  + \beta_t |y| + h(|y|)|z|^{p}  + f(|y|)|z|^2 $ with $0 < p < 2$ and  $f$, $h$  continuous, increasing and positive on $\R_+$.

	\begin{theorem}\label{existenceintegrable} \ 	
		$(i)$ \ Assume that $H$ and $\xi$ satisfy (A1)-(A2).
		Then, $eq(\xi, H)$ has at least one solution $(Y,Z)$ which satisfies for any t,
		\begin{align}\label{bornuy}
		-u_f^{-1}   \left(\E\left [u_{f}\big(\xi_{\alpha, \beta}^-\big)/\mathcal{F}_t \right]\right)
		 \  \leq \
Y_t \
\leq
		 \ u_f^{-1}   \left(\E\left [u_{f}\big(\xi_{\alpha, \beta}^+\big)/\mathcal{F}_t \right]\right).
		\end{align}
 In particular, we have
\begin{align}\label{bornuEy}
		-\E\left [u_{f}\big(\xi_{\alpha, \beta}^-\big)\right]
		 \  \leq \
\E[u_f(Y_t)]
 \ \leq
		 \ \E\left [u_{f}\big(\xi_{\alpha, \beta}^+\big)\right].
\end{align}
		$(ii)$ Among all solutions satisfying \eqref{bornuy}, there are a maximal and a minimal  solution. Note also that, among all solutions satisfying \ $Y^{-g} \leq Y \leq Y^{g}$, there also exists a  maximal and a minimal  solution.
\end{theorem}

	The following remark will be used later. It can be proved as Theorem \ref {existenceintegrable}. 

	\begin{remark} Under the assumptions of Theorem \ref{existenceintegrable}, it also holds that
\begin{align}\label{bornuyalpha}
		-u_f^{-1}   \left(\E\left [u_{f}\big(\xi_{\alpha, \beta}^-\big)/\mathcal{F}_t \right]\right)
\ \leq
		\left(Y_t +\int_0^t\alpha_sds \right)e^{\int_0^t\beta_s ds}  \
		\leq
		 \ u_f^{-1}   \left(\E\left [u_{f}\big(\xi_{\alpha, \beta}^+\big)/\mathcal{F}_t \right]\right).
		\end{align}
		And in particular
\begin{align}\label{bornuEy}
		-\E\left [u_{f}\big(\xi_{\alpha, \beta}^-\big)\right]
 \ \leq \
		 \E\left (u_{f}\left[\left(Y_t + \int_0^t\alpha_sds \right)e^{\int_0^t\beta_s ds}\right] \right)  \
		\leq
		 \ \E\left [u_{f}\big(\xi_{\alpha, \beta}^+\big)\right].
\end{align}

\end{remark}

	\begin{corollary} Theorem \ref{existenceintegrable} remains valid when the function $f$ is continuous but not necessary increasing. In this case, assumption (A2) should be slightly modified as follows:
		
			\vskip 0.2cm\noindent {\rm \textbf{(A2bis)}} \qquad \qquad \qquad \qquad $\displaystyle \E\left[u_{\phi}\big(\xi_{\alpha, \beta}\big)\right] < \infty$ \ \ \ where $\phi(y) := \displaystyle\sup_{0\leq x\leq y} f(x)$.
	\end{corollary}
	\bop Since the function $\phi$ is increasing and continuous, the result follows. \eop
	
	
	\begin{remark}\label{a3relaxed}	Note that when $f$ is merely locally bounded but not necessary continuous, then Theorem \ref{existenceintegrable} remains valid with the following condition in place of assumption (A2).
	
	\vskip 0.2cm\noindent {\rm \textbf{(A2ter)}} \     $\displaystyle \E\left[u_{\varphi}\big(\xi_{\alpha, \beta}\big)\right]  < \infty$,  where $\varphi$  is the smallest continuous,  increasing function such that
	 $\varphi \geq f$.
	
	\end{remark}

\begin{proposition}\label{integrability} (Integrability property of solutions) \ 		
		$(i)$  Let the assumptions of Theorem \ref{existenceintegrable} be satisfied. Assume moreover that there exists $p>1$ such that
		\begin{equation}\label{u(xi)p}
		\E\big[\big(u_{f}[\xi_{\alpha, \beta}]\big)^p\big] < \infty.
		\end{equation}
		then the BSDE$(\xi, H)$ has a solution $(Y, Z)$ which satisfies
		\begin{equation}\label{bornuysup}
		\mathbb{E}\left(\sup_{0\leq t \leq T}\left[u_{f}\bigg(\big[|Y_t| + \int_0^t\alpha_s ds\big]e^{\int_0^t\beta_s ds}\bigg)\right]^p\right)   \  \leq  \
		\mathbb{E}\bigg([u_{f}(\xi_{\alpha, \beta})]^p\bigg).
		\end{equation}

$(ii)$ \  If moreover, $v(|Y|)$ belongs to class $(D)$ and $\sup_{0\leq s \leq T}\E\left[|Y_s|v'(|Y_s|)\right] < \infty$
then, the BSDE$(\xi, H)$ has a solution $(Y, Z)$ which satisfies
\begin{equation}\label{borneZ2}
\E\int_0^T |Z_s|^2 ds < \infty,
\end{equation}
here $v$ is the function defined in Lemma \ref{zvonkinM}-II).
\end{proposition}


\begin{remark} \ $(i)$ \ Let $v$ be the function defined in Lemma
\ref{zvonkinM}-II). Since for $y$ large enough,
  $v'(|y|) \leq \big[u'(|y|)\big]^p$,
$
|y|v'(|y|) \leq \big[u'(|y|)\big]^{p}
$ and we know by Theorem \ref{existenceintegrable} that
\begin{align}\label{bornabsoluey}
|Y_t|\leq
\ u_f^{-1}   \left(\E\left [u_{f}\big(\xi_{\alpha, \beta}\big)/\mathcal{F}_t \right]\right), \nonumber
\end{align}
then the conditions $v(|Y|)$ belongs to class $(D)$ and $\sup_{0\leq s \leq T}\E\left[|Y_s|v'(|Y_s|)\right] < \infty$  are  satisfied when
\begin{equation}\label{u'(xi)palphabeta}
\sup_{0\leq t\leq T}\E\left\{ (\alpha_t + \beta_t)\left( u'_f \left[ u_f^{-1} \left( \E\left [u_{f}\big(\xi_{\alpha, \beta}\big)/\mathcal{F}_t \right]\right)\right]
\right)^p\right\} < \infty,
\end{equation}
and
\begin{equation}\label{u'(xi)p}
\E\left(\big[\sup_{t\leq T}\big( u'_f \big[ u_f^{-1} \left( \E\left [u_{f}\big(\xi_{\alpha, \beta}\big)/\mathcal{F}_t \right]\right)\big]
\big)\big]^{p} \right) < \infty.
\end{equation}

$(ii)$ \ Conditions \eqref{u'(xi)palphabeta} and \eqref{u'(xi)p}, which seem be complicated, cover  those used in  previous works.
For instance, when $\alpha_t$, $\beta_t$ and $f$ are constant with $f(y) = \frac\gamma2$ (see \cite{BH1}), one can take $u_f(y) = exp(\gamma y)$. And in this case, conditions \eqref{u'(xi)palphabeta} and  \eqref{u'(xi)p} become   $\E[e^{p\gamma(\xi_{\alpha, \beta})}] < \infty$ for some $p>1$,  which is  the condition imposed in \cite{BH1}.
\end{remark}



\vskip 0.2cm\noindent	\textbf{Proof of Theorem \ref{existenceintegrable}}.    The existence of (a maximal and a minimal) solutions follows from Proposition \ref{get-gdef} and Lemma \ref{d}. Indeed, the two solutions we constructed in Proposition  \ref{get-gdef} satisfy \ $Y^{-g} \leq  Y^{g}$.   We then use Lemma \ref{d}  with $\xi_1 = -\xi^-$, $\xi_2 = \xi^+$, $H_1 = -g$ and $H_2 = g$ to get the existence of solutions. We shall prove estimate \eqref{bornuy}. For simplicity, we assume that   $\alpha$ and $\beta$ are constant. Since $u_f[e^{\beta T}(\xi^+ + \alpha T)]$ is integrable, then the solution we constructed in the proof of assertion $(i)$  of Proposition \ref{get-gdef} satisfies for any $t$,
	\begin{align*}
	u_f([Y^g +\alpha t]e^{\beta t}) \ \leq  \ \E(u_f[e^{\beta T}(\xi^+ + \alpha T)]/\mathcal{F}_t).
	\end{align*}
	This shows the first inequality of \eqref{bornuy}. A similar argument allows us to prove the second inequality. Theorem \ref{existenceintegrable} is proved. \eop

\vskip 0.2cm \noindent \textbf{Proof of Proposition \ref{integrability}}.
Assertion $(i)$.  Note that the solution $(Y^{5}, Z^{5})$ we constructed in the proof of Proposition \ref{existenceintegrable} satisfies the BSDE	
	\begin{equation}\label{Y5}
	Y_{t}^5= u_f[e^{\beta T}(\xi^+ + \alpha T)] - \int_{t}^{T}Z^5_{s}dW_{s}
	\end{equation} 	
	Since there exists $p>1$ such that inequality
	\eqref{u(xi)p} holds, then  equation \eqref{Y5} has a unique solution $({Y}^5, {Z}^5)$ such that (details can be found in \cite{beh2010, beh2015}):
	
	\begin{equation}\label{bornsupuY5}
	\E(\sup_{0\leq t \leq T} {[Y_t^5]}^p) \leq \E\left(\left[u_f\big(e^{\beta T}[|\xi| + \alpha T]\big)\right]^p\right).
	\end{equation}
But, since $u_f[e^{\beta t}(Y^{g}_t + \alpha_t)] \leq Y^5$, we then have, $\mathbb{E}\big(\sup_{0\leq t \leq T}\big[u_{f}\big(\big[Y_t^g + \int_0^t\alpha_s ds\big]e^{\int_0^t\beta_s ds}\big)\big]^p\big) \leq$

\noindent $\E\left(\left[u_f\big(e^{\beta T}[|\xi| + \alpha T]\big)\right]^p\right)$.
	
	\noindent Similarly, we get
 $\mathbb{E}\big(\sup_{0\leq t \leq T}\big[u_{f}\big(\big[(-Y_t^g) + \int_0^t\alpha_s ds\big]e^{\int_0^t\beta_s ds}\big)\big]^p\big) \leq \E\left(\left[u_f\big(e^{\beta T}[|\xi| + \alpha T]\big)\right]^p\right)$. Assertion (i) is proved.

We shall prove assertion $(ii)$.
Let $v$ be the function defined in Lemma \ref{zvonkinM}-II).
For $N>0$, let
$\tau _{N}:=\inf \{t>0:|Y_{t}|+\int_{0}^{t}|v^{\prime
}(Y_{s})|^{2}|Z_{s}|^{2}ds\geq N\}\wedge T$. Set $\mbox{sgn}(x)=1$ if $x\geq
0 $ and $\mbox{sgn}(x)=-1$ if $x<0$. Since the map $x\mapsto v(|x|)$ belongs to $\mathcal{C}^{2}(
\mathbb{R})$, then thanks to Itô's formula, we have for any $t\in [0,T]$
\begin{align*}
v(|Y_{0}|) & = v(|Y_{t\wedge \tau _{N}}|)
+ \int_{0}^{t\wedge \tau _{N}}\!\!\left[
\mbox{sgn}(Y_{s})v^{\prime }(|Y_{s}|)H(s,Y_{s},Z_{s})-\frac{1}{2}v^{\prime
	\prime }(|Y_{s}|)|Z_{s}|^{2}\right]\!\!ds
\\
& \quad \ \ - \int_{0}^{t\wedge \tau _{N}}\mbox{sgn}(Y_{s})v^{\prime
}(|Y_{s}|)Z_{s}dW_{s} \, .
\end{align*}
Assumption \textbf{(A1)} and
Lemma \ref{zvonkinM}-II) allow us to  show that for any $N>0$,
\begin{align}\label{estimz}
\frac{1}{2}\mathbb{E}\!\int_{0}^{t\wedge \tau _{N}}\!\!|Z_{s}|^{2}ds
&\leq
\E\big[v(|Y_{t\wedge \tau _{N}}|)\big]+\mathbb{E}\!\int_{0}^{T}\!\!\left[\alpha_t  v^{\prime
}(|Y_{s}|)+\beta_t |Y_{s}|v^{\prime
}(|Y_{s}|)\right]\!ds
\end{align}
Since  the processes  $v(|Y|)$ and $|Y|v'(|Y|)$ belong to class $(D)$,
the proof is completed by using Fatou's lemma.
Proposition \ref{integrability} is proved. \eop


\begin{corollary}\label{bmo} (BMO property) $(i)$ \  Let (A1) be satisfied. Assume moreover that \  $\xi$, $\int_0^T\alpha_sds$ and $\int_0^T\beta_sds$ are bounded.  Then, every solution $(Y, Z)$ satisfying  inequalities \eqref{bornuy} is such that $Y$ is bounded and the process  $(\int_0^t Z_sdW_s)_{0\leq t \leq T}$ is a BMO martingale.
	
\end{corollary}


\bop \
$(i)$ \ Let $(Y, Z)$ be a solution to  $eq(\xi, H)$ such that  $Y$  satisfies  inequalities \eqref{bornuy}. Since $\xi$, $\int_0^T\alpha_sds$ and $\int_0^T\beta_sds$ are bounded, then clearly $Y$ is bounded. Arguing as in the proof of Proposition \ref{integrability} $(ii)$, one can show that $Z$ belongs to $\mathcal{M}^2$. We shall prove that the process  $(\int_0^t Z_sdW_s)_{0\leq t \leq T}$ is a BMO martingale.
Let $v$ be
the function defined in Lemma
\ref{zvonkinM}-II). Since the map $x\longmapsto v(|x|)$ belongs to $\mathcal{C}^{2}(
\mathbb{R})$, then  It\^o's formula shows that for
any $\mathcal{F}_t$--stopping time $\tau \leq T$,
\begin{align*}
v(|Y_{ T}|)& = v(|Y_{\tau}|) +  \int_{\tau}^{T}\!\!\left[
\frac{1}{2}v^{\prime \prime }(|Y_{s}|)|Z_{s}|^{2}-\mbox{sgn}(Y_{s})v^{\prime
}(|Y_{s}|)H(s,Y_{s},Z_{s})\right]\!\!ds \\
& \ \ \ \ \ + \int_{\tau}^{T}\mbox{sgn}(Y_{s})v^{\prime
}(|Y_{s}|)Z_{s}dW_s \, .
\end{align*}
Since $Y$ is bounded and $Z$ belongs to $\mathcal{M}^2$, it follows that the stochastic
integral in the right hand side term of the previous equality is a square integrable $\mathcal{F}_t$--martingale.
Passing to conditional expectation, one can show that there exist positive constants $K_1$ and $K_2$  such that
\begin{align*}
\mathbb{E}\left(\int_{\tau}^{T}|Z_{s}|^{2}ds/\mathcal{F}_{\tau}\right)
& \leq
K_1 + \mathbb{E}\left(\!\int_{\tau}^{T}\!\!\left[(\alpha_s+\beta_s|Y_{s}|)v^{\prime
}(|Y_{s}|)\right]\!ds/\mathcal{F}_{\tau}\right)
\end{align*}
We complete the proof by noticing that the processes $Y$,  $\int_0^T\alpha_s ds$ and $\int_0^T\beta_s ds$   are bounded.   Assertions $(ii)$ and $(iii)$ can be proved similarly. \eop


\begin{remark} Let $|H(t,y,z)| \leq \alpha_t + \beta_t |y| + f(y)|z|^2$, with $f$ continuous, positive and increasing but not globally integrable. Then, the condition which ensures the existence of solution is:

	\begin{equation}\label{f(0)}
	\E\big(u_f(\xi_{\alpha, \beta}^+)) < \infty \ \ \ \hbox{and} \ \ \  \E\bigg(\exp\bigg[(\xi^-)e^{\int_0^T\beta_sds}2f(0)\bigg]\bigg) < \infty
	\end{equation}
	
\end{remark}

\bop \ We assume that $\alpha$ and $\beta$ are constants for simplicity. As previously, thanks to  Lemma \ref{d}, $eq(\xi, H)$ has a solution if $eq(\xi^+, \alpha + \beta |y| + f(y)|z|^2)$ has a positive solution and $eq(-\xi^-, -[\alpha + \beta |y| + f(y)|z|^2])$ has a negative solution. The first inequality in \eqref{f(0)} can be proved as in Proposition  \ref{get-gdef}. To prove the second inequality of \eqref{f(0)}, we  see that $(Y,Z)$ is a negative solution to $eq(-\xi^-, -[\alpha + \beta |y| + f(y)|z|^2])$   if  $(Y', Z') := (-Y,-Z)$ is a positive solution to the BSDE 	
\begin{equation*}
Y_t' = \xi^- + \int_t^T [\alpha + \beta Y_s' + f(-Y_s')|Z_s'|^2]ds - \int_t^T Z_s'dW_s.
\end{equation*}

\noindent Since $f$ is increasing, we use  Lemma \ref{d} to show that the previous BSDE has a positive solution if the BSDE $Y^{1}_t = \xi^- + \int_t^T [\alpha + \beta Y^{1}_s + 2f(0)|Z^{1}_s|^2]ds - \int_t^T Z^{1}_s dW_s$ has a positive solution. The sequel of the proof goes as in Theorem \ref{existenceintegrable}.
\eop


\vskip 0.2cm The following corollary extends the result of \cite{beo2013, beo2017} to the case where $\xi_{\alpha, \beta}$ is only integrable and the coefficients $\alpha$, $\beta$, $\gamma$ are positive processes. 
\begin{corollary}\label{aop} 
	Let (A1) be satisfied with $f$ globally integrable and locally bounded. Assume that $\xi_{\alpha, \beta}$ is integrable. Then, $eq(\xi, H)$ has a solution such that for every t 
	\begin{align}\label{bornyaop}
	- \E\left (\xi_{\alpha, \beta}^-/\mathcal{F}_t \right)
	\  \leq \
	Y_t \
	\leq
	\E\left (\xi_{\alpha, \beta}^+/\mathcal{F}_t \right).
	\end{align}
	$(ii)$ Among all solutions satisfying \eqref{bornuy}, there are a maximal and a minimal  solution. Note also that, among all solutions satisfying \ $Y^{-g} \leq Y \leq Y^{g}$, there also exists a  maximal and a minimal  solution where the function $g$ is defined by \eqref{g}.	
\end{corollary}

\bop
Let $g$ be the function defined by \eqref{g}. 

 Since $f$ is globally integrable $u_f$ and its inverse are uniformly Lipschitz. Hence,  arguing as in the proof of Proposition \ref{get-gdef}  one show that $eq(\xi^+, g)$ has a positive solution  $Y^{g}$ such that:  
$0\leq Y_t^{g} \
\leq
\E\left (\xi_{\alpha, \beta}^+/\mathcal{F}_t \right)
$. Symmetrically, we show that $eq(-\xi^-, -g)$ has a negative solution  $Y^{-g}$ which satisfied : 
$
- \E\left (\xi_{\alpha, \beta}^-/\mathcal{F}_t \right)
\  \leq \
Y_t^{-g} \ \leq 0$. 
  Since $f$ is locally bounded, then according to the proof  of Theorem 4.1 of  \cite{beo2017}, there exists a continuous function $\bar {f}$ such that $f\leq \bar {f}$. Therefore, we can apply  Lemma \ref{d} with \ $\eta_t = \alpha_t + \beta_t(|Y_t^{-g}| + |Y_t^{g}|) $
  and
  $C_t =  
  \sup_{s\leq t}\sup_{a\in [0,1]}|\bar{f}(a 
  Y_s^{-g} +(1-a) Y_s^{g})|$. Corollary \ref{aop} is proved.   \eop



\subsubsection{ BSDE$(\xi, H)$ with $|H(t,y,z)| \leq \alpha_t + \beta_t |y| + \theta_t |z| + f(|y|)|z|^2$}



As in the previous subsection, we use Lemma \ref{d} to reduce  the solvability of $eq(\xi, H)$  to the positive solvability ($Y \geq 0$) of the simple BSDEs $(u_f(\xi_{\alpha, \beta}^+), \theta_t|z|)$ and $(u_f(\xi_{\alpha, \beta}^-), \theta_t|z|)$ then apply Proposition  \ref{propehz} to conlude.  We put 

\begin{equation}\label{hh} 
h(t,y,z):= \alpha_t + \beta_t |y| + \theta_t |z| + f(|y|)|z|^2
\end{equation}

 

\begin{proposition}\label{propz}
Assume that {\rm{(A3)}}, {\rm{(A4)}} are satisfied. Then,

$(i)$ \ $eq(\xi, H)$ has at least one solution such that
\begin{equation}\label{bornhz}
 - u_f^{-1}\left(\displaystyle{ess\sup_{\pi\in \sum}} \E\left(\Gamma_{t,
	T}^{\pi}{u_f(\xi_{\alpha, \beta}^-)}/{\cal F}_t\right)\right)
\leq 	
Y_t
\leq  u_f^{-1}\left(\displaystyle{ess\sup_{\pi\in \sum}} \E\left(\Gamma_{t,
		T}^{\pi}{u_f(\xi_{\alpha, \beta}^+)}/{\cal F}_t\right)\right)
\end{equation}

$(ii)$ \ Among all solutions satisfying inequalities \eqref{bornhz} there are a maximal and a minimal solution.

$(iii)$ \ Assume moreover that $\exp(\int_0^T\theta_s^2ds)$ is integrable,  and $\xi_{\alpha, \beta}$ is bounded. Then all solutions satisfying inequalities \eqref{bornhz} are bounded. 
\end{proposition} 

\begin{remark}
	 One may wonder what is the usefulness of assumption (A3) since it can be reduced to assumption (A1) by the operation $\alpha_t + \beta_t |y| + \theta_t |z| + f(|y|)|z|^2)\leq \alpha_t + \frac12 \theta^2 + \beta_t |y| +   [\frac12 + f(|y|)]|z|^2)$. 
	  It should be noted that in this case the integrability requested to the terminal value will be higher.
\end{remark}


\vskip 0.2cm \noindent\textbf{Proof of Proposition \ref{propz}. }
Since $H$ satisfies (A3), then according to  Lemma \ref{d}, it is enough to show that $eq(\xi^+, \alpha_t + \beta_t |y| + \theta_t |z| + f(|y|)|z|^2)$ has a positive solution.   Arguing as in the proof of Theorem  \ref{existenceintegrable}, one can show that $eq(\xi^+, \alpha_t + \beta_t |y| + \theta_t |z| + f(|y|)|z|^2)$ has a positive solution if $eq(u_f(\xi^+_{\alpha, \beta}),\,   \theta_t |z|)$ has a  solution which is greater than $u_f(\alpha T e^{\beta T})$.  This implies, thanks to  Proposition \ref{propehz}, that  assumption (A4) is satisfied. Proposition \ref{propz} is proved.  \eop


\begin{remark}  
 $(i)$ \ 	Taking $f=0$ in  Proposition \ref{bornhz}, we get the existence of one dimensional BSDEs with a stochastic linear growth. This covers the results of \cite{LSM1} and \cite{ht}.  
 
 $(ii)$ We emphasize that Proposition \ref{ehz} combined with Lemma \ref{d} allows to  directly prove the existence of solutions to BSDEs whose generators satisfy \begin{equation}\label{Hlinz} 
 |H(t,y,z)| \leq \alpha_t + \beta_t |y| + \theta_t |z| 
 \end{equation} 
 and $\xi$ satisfies assumption (A3) with $f=0$. 
 
 This constitute a new result on the existence of solutions to BSDEs with at most linear growth which for instance covers the recent result \cite{ht}, with a simpler  proof.  
\end{remark}

\begin{remark}\label{existavecz}
 Since the transformation $u_f$ does not impact $Z$, we then have
	
 \ If \ $|H(t, y, z)| \leq \alpha_t + \beta_t|y| + \theta_t|z| + f(|y|)|z|^2$ with $\alpha$, $\beta$, $f$ satisfy (A3) and $\E\int_0^T e^{q\gamma_s}ds < \infty$ for some $q>0$, then arguing as in Proposition \ref{propz} and using   \cite{beh2010, beh2015}, one can show that $eq(\xi, H)$ has a solution $(Y,Z)$ in $\mathcal{S}^p\times\mathcal{M}^p$ provided that
$\exp(\frac12\int_0^T\gamma_se^{\beta_s}ds)u_f\left(\xi_{\alpha, \beta}\right)$ is $p$--integrable for some $p>1$.
\end{remark}


\begin{corollary}\label{bmoz} (BMO property) $(i)$ \  Let (A3) be satisfied. Assume moreover that \  $\xi$, $\int_0^T\alpha_sds$, $\int_0^T\beta_sds$ and $\int_0^T\gamma_s^2ds$  are bounded.  Then, every solution $(Y, Z)$ satisfying  inequalities \eqref{bornhz} is such that $Y$ is bounded and the process  $(\int_0^t Z_sdW_s)_{0\leq t \leq T}$ is a BMO martingale.
	
 $(ii)$ \ When $H$ is dominated by $\alpha_t + \theta_t|z| + f(|y|)|z|^2$ with $\xi$, $\int_0^T\alpha_sds$ and $\int_0^T\gamma_s^2ds$  bounded, then we have the same conclusion as $(i)$ with $f$  locally integrable and increasing but not continuous.

$(iii)$ \ When $H$ is dominated by $\theta_t|z| + f(|y|)|z|^2$ with $\xi$, $\int_0^T\alpha_sds$ and $\int_0^T\gamma_s^2ds$  bounded, then we have the same conclusion as $(i)$ with $f$  merely locally integrable but neither increasing nor continuous.

\end{corollary}

	
	\vskip 0.2cm The following Corollary can be proved by combining the proof of Corollary \ref{aop}  with that of Proposition \ref{propz}. 
	
	\begin{corollary}\label{aopz} 
		Let (A3) be satisfied with $f$ globally integrable and locally bounded. Assume that $\xi_{\alpha, \beta}$ satisfies 
		
		$$
		 \qquad \displaystyle  \sup_{\pi\in\sum} \E\left(\Gamma_{0, T}^{\pi}(\xi_{\alpha, \beta})\right) :=\sup_{\pi\in\sum} \E\left({}e^{\int_0^T \theta_u \pi_udW_u -\frac12 \int_0^T
			\theta_u^2 |\pi_u|^2du}(\xi_{\alpha, \beta})\right)  < +\infty
		$$ 
		 
		Then, $eq(\xi, H)$ has a solution such that for every t 
	\begin{equation}\label{bornhzaop}
	- \displaystyle{ess\sup_{\pi\in \sum}} \E\left(\Gamma_{t,
		T}^{\pi}{(\xi_{\alpha, \beta}^-)}/{\cal F}_t\right)
	\leq 	
	Y_t
	\leq  \displaystyle{ess\sup_{\pi\in \sum}} \E\left(\Gamma_{t,
		T}^{\pi}{(\xi_{\alpha, \beta}^+)}/{\cal F}_t\right)
	\end{equation}
	
	$(ii)$ \ Among all solutions satisfying inequalities \eqref{bornhzaop} there are a maximal and a minimal solution. 
	
		 Note also that, among all solutions satisfying \ $Y^{-h} \leq Y \leq Y^{h}$, there also exists a  maximal and a minimal  solution where the function $g$ is defined by \eqref{hh}.	
	\end{corollary}




\section { Quadratic BSDEs and BSDEs with logarithmic nonlinearity}
The aim of this subsection is to study the BSDE $(\xi, H)$ with $H$ continuous and with $lLogl$-growth. That is : there exist positive constants $a$, $b$ and $c$ such that for every $t, y, z$
\begin{equation}\label{Hlog}
|H(t, y, z) \leq a + b |y| + c |y||\ln |y||
\end{equation}
Using  Lemma \ref{d}, we show that that of $eq(\xi, H)$ is equivalent to the solvability of $eq(\xi, \alpha + \beta\vert y\vert + \frac{\gamma}{2}|z|^2)$, for suitable $\alpha$, $\beta$ and $\gamma$. According to Theorem \ref{existenceintegrable},  the BSDE $(\xi, \alpha + \beta\vert y\vert + \frac{\gamma}{2}|z|^2)$ has a solution when $exp(\gamma e^{\alpha T}|\xi|)$ is integrable, and we have the following propositions.

\begin{proposition}\label{QBSDEZ2log}
	The BSDE $(\xi, \alpha + \beta\vert y\vert + \frac{\gamma}{2}|z|^2)$ has a solution if and only if the BSDE $(e^{\gamma\xi}, \gamma\alpha y + \gamma\beta y|\ln y|)$ has a solution.
\end{proposition}

\bop \   Let
\begin{equation}\label{glog}
g(t,y,z):= \alpha + \beta\vert y\vert + \frac{\gamma}{2}|z|^2.
\end{equation}
Applying Itô's formula to $u(y) := e^{\gamma y}$, we  show that $(Y_t, Z_t)$ is a solution to $eq(\xi, g)$ if and only if $(\bar{Y_t}, \bar {Z_t}) := (e^{\gamma Y_t}, \gamma e^{\gamma Y_t}Z_t) $ is  a solution to $eq(e^{\gamma\xi}, \gamma\alpha y + \gamma\beta y|\ln y|)$. Indeed:
Let  $(Y, Z)$ be
a  solution of   $eq(\xi, g)$. By Itô's formula we have,
\begin{eqnarray*}
	e^{\;\gamma Y_t}&=&e^{\;\gamma Y_T}+\int_t^T \gamma e^{\gamma Y_s}g(s,Y_s,Z_s)ds
	-\int_t^T \gamma  e^{\gamma Y_s}Z_sdW_s-\frac{\gamma^2}{2}\int_t^Te^{\gamma Y_s}|Z_s|^2ds
	\\
	&=&e^{\,\gamma \xi}+\int_t^T \gamma e^{\gamma Y_s}(\alpha +\beta |Y_s|+\frac{\gamma}{2}|Z_s|^2)ds
	- \int_t^T\gamma  e^{\gamma Y_s}Z_sdW_s-\frac{\gamma^2}{2}\int_t^Te^{\gamma Y_s}|Z_s|^2ds\\
	&=&e^{\,\gamma\xi}+\int_t^T  \gamma e^{\gamma Y_s}(\alpha + \beta|Y_s|)ds-\int_t^T \gamma e^{\gamma Y_s}Z_sdW_s
\end{eqnarray*}
\\
It is clear that $\bar Y$ $>$
$0$ and  $(\bar Y, \bar Z)$ satisfies the BSDE
\begin{equation*}\label{eqbar}
\bar{Y_t}=e^{\;\gamma\xi}+\int_t^T  \gamma \left(\alpha\bar{Y_s}+\beta\bar{Y_s}|\ln\bar
Y_s|\right)ds-\int_t^T\bar{Z_s}dW_s.
\end{equation*}
Proposition \ref{QBSDEZ2log} is proved.
\eop


\begin{proposition}\label{ylogyG}
	Let  $\xi$ be an $\mathcal{F}_T-$measurable random variable.  Let $G$ be defined by
	\begin{equation}\label{G}
	G(y): = a + b\vert y\vert  + c|y||\ln|y||
	\end{equation}
	
	$(i)$ \  If \  $eq(e^{(a + b  + 2c)  e^{c\, T}}(\xi^+ + 1)^{ e^{c\, T}},\ 0)$ has a positive solution, then \ $eq(\xi^+, G)$ has a positive solution.
	
	$(ii)$ \ $(Y, Z)$ is a negative solution of $eq(-\xi^-, -G)$  if and only if $(-Y, -Z)$ is a positive solution to $eq(\xi^-, G)$. Therefore,  if  \  $eq(e^{(a + b  + 2c)  e^{c\, T}}(\xi^- + 1)^{ e^{c\, T}},\ 0)$ has a positive solution, then \ $eq(-\xi^-, -G)$ has a negative solution.
		
		$(iii)$ \ If $|\xi|^{e^{cT}}$ is integrable, then $eq(\xi^+, G)$ has a positive solution and $eq(-\xi^-, -G)$ has a negative solution, and  therefore $eq(\xi, H)$ has at least one solution $(Y,Z)$ which belongs to $\mathcal{S}^{e^{cT}}\times\mathcal{M}^2$. Moreover, according to see \cite{bcimpa},  the uniqueness holds in $\mathcal{S}^{e^{2cT} + 1}\times\mathcal{M}^2$ provided that $|\xi|^{e^{2cT} + 1}$ is integrable.
\end{proposition}

\bop
Let $G(y): = a + b\vert y\vert  + c|y||\ln|y||$.  Let $Y^G$ be a positive solution to $eq(\xi^+, G)$. This is equivalent to say that $Y^G$ is a positive solution to the BSDE $(a + b y  + cy |\ln y|)$. Applying Itô's formula to the function $u(Y^G_t) := \ln (Y^G_t+1)$, we obtain
\begin{align}\label{ln(y+1)}
u(Y^G_t) =  \ln(\xi^+ + 1) & + \int_t^T \left( \big[ a + b Y^G_s + c\,Y^G_s|\ln(Y^G_s)|\big]\frac{1}{1 + Y^G_s} +  \frac{1}{2}\frac{|Z^G_s|^2}{(1 + Y^G_s)^2}\right)\,ds \notag
\\
& - \int_t^T \frac{1}{1 + Y^G_s}\, Z^G_s dW_s.
\end{align}
The process $(\bar{Y}, \bar{Z}) := (\ln(1 + Y^G), \frac{Z^G}{1 + Y^G})$ satisfies the BSDE
\begin{align}\label{log Y+1}
\bar{Y}_t =  \ln(\xi^+ + 1)  + \int_t^T  \bar{H}(s, \bar{Y}_s, \bar{Z}_s) ds
- \int_t^T \bar{Z}_s\,dW_s \, ,
\end{align}
where
\begin{align}\label{Hbar}\displaystyle \bar{H}(t,y,z) := \big[ a + b (e^{y}-1) + c\,(e^{y}-1)\,|\ln(e^{y}-1)|\big]\frac{1}{e^{y}} +  \frac{1}{2}|z|^2.
\end{align}
Since the function $x|\ln(x)| < 1$ for each $x$ in $[0,\, 1]$ and  strictly increasing in  $[1, \ +\infty)$, we then have $
\big(x |\ln(x)|\big)\frac{1}{1 + x} \leq  1 + |\ln(x+1)|
$. Hence
\begin{align}\label{boundyG}
\big(a + b x + c x |\ln(x)|\big)\frac{1}{1 + x} \leq a + b + c + c|\ln(x+1)|
\end{align}
It follows that
\begin{align}\label{leqHbar}
0\, \leq\, \bar{H}(t,y,z)\, \leq \,   a + b\,  + c + c\,y +  \frac{1}{2}|z|^2
\end{align}
According to Lemma \ref{d}, it is enough to show that $ eq(\ln(\xi^+ + 1), a + b  + c + c\,y +  \frac{1}{2}|z|^2)$ has a positive solution.
But from Proposition \ref{get-gdef},  \ $eq(\ln(\xi^+ + 1), a + b  + c + c\,y +  \frac{1}{2}|z|^2)$ has a positive solution when \ $eq(e^{(a + b  + 2c)  e^{c\, T}}e^{ e^{c\, T}\ln(\xi^+ + 1)},\ 0)$ \  has a positive solution, which is equivalent to say that $eq(e^{(a + b  + 2c)  e^{c\, T}}(\xi^+ + 1)^{ e^{c\, T}},\ 0)$ has a positive solution. This implies  that $(\xi^+ + 1)^{ e^{c\, T}}$ is integrable. Assertions $(ii)$ can be proved as assertion $(i)$. Lemma \ref{d} allows to establish existence of solutions of assertion $(iii)$.
\eop


\begin{remark}\label{uniquezcarre} {\rm{(Uniqueness)}}.
According to assertion $(iii)$ of Proposition \ref{ylogyG}, $eq(\xi, a + b |y| + c |y||\ln |y||)$ has a unique solution $(Y, Z)$ which belongs to $\mathcal{S}^{e^{2cT} + 1}\times\mathcal{M}^2$ provided that $|\xi|^{e^{2cT} + 1}$ is integrable.  Therefore,   $eq(\xi, \alpha + \beta|y| + \frac{\gamma}{2}|z|^2)$ has a unique solution provided that
 $\exp\big({\;\gamma\xi(e^{2\beta T }+1)}\big)$ is integrable. We moreover have,  $\displaystyle \sup_{0\leq s\leq T}\exp\left(\gamma|Y_s|(e^{2\beta_s}+1)\right)$ is integrable. This gives a simple proof to the uniqueness of  $eq(\xi, \alpha + \beta|y| + \frac{\gamma}{2}|z|^2)$ without using the convexity  (in $z$) of the  generator.
\end{remark}


\begin{remark}  The uniqueness of solutions under assumptions (A1)-(A2) as well as under assumptions (A3)-(A4), and the existence of viscosity solutions  to the related partial differential equation are in progress. 
\end{remark}

%
%

\section{Appendix.}
\vskip 0.1cm   

We recall the result of	Essaky $\&$ Hassani (\cite{EH2013}) on the two barriers reflecting QBSDEs. It establishes the existence of solutions for reflected QBSDEs without assuming any integrability condition on the terminal datum. This result is used in the proof of Lemma \ref{d}. 

\begin{theorem}\label{EH2013}
	{\rm (\cite{EH2013}, Theorem 3.2)} Let $L$ and $U$ be continuous
	processes and $\xi$ be a $\mathcal{F}_T$--measurable random variable. Assume
	that
	
	1) \hskip 0.2cm for every $t\in[0, \ T]$, \ $L_t \leq U_{t}$
	
	2) \hskip 0.2cm $L_T\leq\xi\leq U_T$.
	
	3) \hskip 0.2cm there exists a continuous semimartingale which passes
	between the barriers $L$ and $U$.
	
	4) \hskip 0.2cm The generator $h$ is continuous in $(y,z)$ and satisfies for every $
	(s,\omega)$, every $y\in [L_s(\omega), U_s(\omega)]$ and every $z\in \mathbb{
		R}^d$.
	
	\vskip 0.2cm
	\hskip 3cm $|h(s, \omega, y, z )| \leq \eta_s(\omega)+ C_s(\omega)|z|^2 $
	
	\vskip 0.2cm\noindent where $\eta$ and $C$ are two $\mathcal{F}_t$--adapted  processes such that
	$\E\int_0^T \eta_s ds < \infty$ and $C$ is  continuous.
	
	\vskip 0.2cm\noindent Then, the following RBSDE has a maximal and a minimal 
	solution.
	\begin{equation}  \label{eq000}
	\left\{
	\begin{array}{ll}
	(i) & Y_{t}=\xi + \displaystyle\int_{t}^{T}h(s,Y_{s},Z_{s})ds  - \displaystyle
	\int_{t}^{T}Z_{s}dW_{s} \\
	& + \displaystyle
	\int_{t}^{T}dK_{s}^+ - \displaystyle\int_{t}^{T}dK_{s}^- \; \mbox{for all }\; t\leq T \\
	(ii) & \forall \ t\leq T,\,\, L_t \leq Y_{t}\leq U_{t},\quad \\
	(iii) & \displaystyle\int_{0}^{T}( Y_{s}-L_{s}) dK_{s}^+ = \displaystyle
	\int_{0}^{T}( U_{s}-Y_{s}) dK_{s}^-=0,\,\, \mbox{a.s.}, \\
	(iv) & K_0^+ =K_0^- =0, \,\,\,\, K^+, K^- \,\,
	\mbox{are continuous
		nondecreasing}. \\
	(v) & dK^+ \bot dK^-
	\end{array}
	\right.
	\end{equation}
\end{theorem}


\vskip 0.3cm The following lemma allows to remove the quadratic term from $eq(\xi, \alpha_t + \beta_t|y| + \theta_t|z| + f(|y|)|z|^2$

\begin{lemma}
	\label{zvonkinM}  {\rm I)} \ Let $f \in \mathbb{L}^{1}_{loc}(
	\mathbb{R})$ but not necessarily continuous. Then the
	function
	\begin{equation}\label{uzvonkin0}
	u_f(x):=\int_{0}^{x}\exp \left( 2\int_{0}^{z}f(r)dr\right) dz
	\end{equation}
	satisfies  the differential equation $
	\frac{1}{2}u_f^{\prime \prime }(x)-f(x)u_f^{\prime }(x)=0$ \ $a.e$  on $\mathbb{R}$, and has the following properties:
	
	(j) \ $u_f$ is a one to one function. Both $u$ and its inverse $u_f^{-1}$ are locally Lipschitz, that is for every $R>0$ there exist two positive constants $
	m_R $ and $M_R$ such that, for any \ $|x|, |y| \leq R$,
	
	\hskip 1.5cm  \ $m_R\left\vert x-y\right\vert
	\leq \left\vert u_f(x)-u_f(y)\right\vert \leq M_R\left\vert x-y\right\vert $
	
	(jj)  \  Both \ $u_f$ and its inverse function $u_f^{-1}$ belong to $ W_{1,\, loc}^{2}(\mathbb{R
	})$. If moreover $f$ is continuous, then both $u_f$ and $u_f^{-1}$ belong to $\mathcal{C}^{2}(\mathbb{R})$.

	\vskip 0.2cm {\rm II)} \ Set
	$$K(y):=\int _{0}^{y} \exp\left(-2\int_{0}^{z}f(r)dr\right)dz.$$
	Then, the function \begin{equation}\label{vzvonkin1}
	v(x):=\int_{0}^{x}K(y) \exp\left(2\int_{0}^{y}f(r)dr\right)dy
	\end{equation}
	satisfies the differential equation \ $\frac{1}{2}
	v^{\prime\prime}(x)- f(x)v^{\prime}(x)= \frac{1}{2}$ \ $a.e.$ on $\R$ and has the following properties:
	
	\vskip 0.15cm  $(jjj)$   $v$ and $v^{\prime }$ are positive on $\R_+$  \ and \ $v$ belongs to $W_{1,\, loc}^{2}(\mathbb{R})$.
	
	\vskip 0.15cm \noindent
	
	\vskip 0.15cm  $(jv)$ The map $x\longmapsto  v(\vert x\vert)$ belongs
	to $W _{1,\, loc}^{2}(\mathbb{R})$, and belongs to $\mathcal{C}^{2}(\mathbb{R})$ if $f$ is continuous.
	
	
	\vskip 0.2cm {\rm III)} \ Set \ $G(z) := \int_0^z f(x)e^{-2\int_0^x f(r)dr}dx$. The function
	\begin{equation}
	\label{wzvonkinf}
	w(y) :=  \int_0^y G(z)e^{2\int_0^x f(r)dr}dz
	\end{equation}
	has the following properties :
	
	\vskip 0.15cm(vj) \  the map $x\longmapsto w(|x|)$ belongs to $W _{1,\,loc}^{2}$
	
	\vskip 0.15cm (vjj) \ $w$ satisfies the following differential equation
	\begin{equation}\label{w}
	\frac12 w^{\prime\prime}(x) - f(x) w'(x) = \frac12 f(x) \quad  a. e. \hskip 0.2cm on \hskip 0.1cm \R.
	\end{equation}
\end{lemma}

\bop
\ I) \  Clearly, $u_f$ and its inverse $u_f^{-1}$ are continuous, one to one,  strictly increasing
functions and we have \ $
\frac{1}{2}u_f^{\prime \prime }(x)-f(x)u_f^{\prime }(x)=0$ \ $a.e.$ on $\R$.
 Since \
$u_f'(x):= \exp (2\int_{0}^{x}f(t)dt)$,  then, 
\begin{equation}\label{quasiisometrylocale}
\mbox{for every} \ |x|\leq R, \ \
\exp\left(-2\| f\|_{\L^1([-R,\ R])}\right)\ \leq \ |u_f'(x)| \ \leq \ \exp\left(2\| f\|_{\L^1([-R,\ R])}\right).
\end{equation}
This shows that  $u_f$ and $u_f^{-1}$ are locally Lipschitz.

We prove $(jj)$. Using inequality \eqref{quasiisometrylocale}, one can show that both $u_f$
and $u_f^{-1}$ belong to $\mathcal{C}^1$. Since the second generalized derivative $u_f^{\prime \prime
}$ satisfies $u_f^{\prime \prime
}(x) = 2f(x)u_f^{\prime}(x)$ for $a.e. \ x$,  we get that $u_f^{\prime \prime
}$ belongs to ${\L}_{loc}^1(\R)$. Therefore $u_f$ belongs to  $W_{1,\, loc}^{2}(\mathbb{R})$.
Using again assertion $(j)$, we prove that $u_f^{-1}$ belongs to $W_{1,\, loc}^{2}(\mathbb{R})$.

II) Obviously $v$ and $v^{\prime }$ are positive on $\R_+$ and $v$ satisfies the differential equation
\ $\frac{1}{2}
v^{\prime\prime}(x)- f(x)v^{\prime}(x)= \frac{1}{2}$ \ $a.e.$ on $\R$.
Since $f$ is locally integrable on $\R$, one can easily  check that
$v$ belongs to $W_{1,\, loc}^{2}(\mathbb{R})$. This proves assertions $(jjj)$, from which we deduce  assertion $(jv)$. The proof of III) is similar.
\eop


\vskip 1cm\noindent 
\textbf{Acknowledgment.} The author would like to express special thanks to Rainer Buckdahn, El-Hassan Essaky, Said Hamadene, Mohammed Hassani and Ludovic Tangpi for various discussions about this work.

 This work has been presented at Conference on Stochastic control, BSDEs and new developments (September 2017 at Roscoff, France), International Conference on Stochastic Analysis (October 2017 at Hammamet, Tunisia), Journées Analyse Stochastique (January 2018, Biskra, Algeria), Conference SMT (March 2018, Tabarka, Tunisia), and Séminaire Bachelier (January 2018, IHP, Paris, France),  Séminaire de l'université du Mans (February 2018, Le Mans, France), Viennese Seminar on Probability Theory and Mathematical Finance (Vienna, July 2018) and Workshop on Probability (September 2018, Marrakech, Maroc).  The author is sincerely  grateful  to the organizers of these events for their  invitation.  


\end{document}